\documentclass[pdflatex,sn-mathphys-num]{sn-jnl}% Math and Physical Sciences Numbered Reference Style

\usepackage{graphicx}%
\usepackage{multirow}%
\usepackage{amsmath,amssymb,amsfonts}%
\usepackage{amsthm}%
\usepackage{mathrsfs}%
\usepackage[title]{appendix}%
\usepackage{xcolor}%
\usepackage{textcomp}%
\usepackage{manyfoot}%
\usepackage{booktabs}%
\usepackage{algorithm}%
\usepackage{algorithmicx}%
\usepackage{algpseudocode}%
\usepackage{listings}%
\usepackage{bm}
\usepackage{subcaption}
\usepackage{setspace}
\usepackage{pgfplots}
\pgfplotsset{compat=1.17}
\usetikzlibrary{external}
\usepgfplotslibrary{external}

\theoremstyle{thmstyleone}%
\theoremstyle{thmstyletwo}%

\theoremstyle{thmstylethree}%

\def\ra{$r$-adaptivity}

\def\bx{{\boldsymbol{x}}}   %bold x
\def\bxi{{\boldsymbol{x}_{\tt init}}}

\def\lim{{\xi}}    %limiting function

\def\bbx{\bar{\bx}}
\def\br{\mathbf{r}}   %bold r

\def\mcm{\mathcal{M}}

\def\mcf{F}
\def\mcr{R}
\def\mcu{\mathcal{U}}

\definecolor{code}{rgb}{0.7, 0, 0.4}

\def\bmiv{\boldsymbol{\mathit{v}}}
\def\bb{{\boldsymbol{b}}}   %bold b
\def\bu{{\boldsymbol{\mathit{u}}}}   %bold u
\def\bsig{{\boldsymbol{\sigma}}}
\def\bw{\boldsymbol{w}}
\def\tbw{\tilde{\boldsymbol{w}}}

\begin{document}

\title[PDE-Constrained High-Order Mesh Optimization]{PDE-Constrained High-Order Mesh Optimization}

%%  \sfx{IV}}\email{iauthor@llnl.gov}
%%=============================================================%%

\author[1]{\fnm{Tzanio} \sur{Kolev}}\email{kolev1@llnl.gov}
\author[1]{\fnm{Boyan} \sur{Lazarov}}\email{lazarov2@llnl.gov}
\author[1]{\fnm{Ketan} \sur{Mittal}}\email{mittal3@llnl.gov}
\author*[1]{\fnm{Mathias} \sur{Schmidt}}\email{schmidt43@llnl.gov}
\author[1]{\fnm{Vladimir} \sur{Tomov}}\email{tomov2@llnl.gov}
\affil[1]{\orgname{Lawrence Livermore National Laboratory}, \orgaddress{\street{7000 East Avenue}, \city{Livermore}, \postcode{94550}, \state{California}, \country{USA}}}


\abstract{We present a novel framework for PDE-constrained $r$-adaptivity of high-order meshes. The proposed method formulates mesh movement as an optimization problem, with an objective function defined as a convex combination of a mesh quality metric and a measure of the accuracy of the PDE solution obtained via finite element discretization. The proposed formulation achieves optimized, well-defined high-order meshes by integrating mesh quality control, PDE solution accuracy, and robust gradient regularization. We adopt the Target-Matrix Optimization Paradigm to control geometric properties across the mesh, independent of the PDE of interest. To incorporate the accuracy of the PDE solution, we introduce error measures that control the finite element discretization error. The implicit dependence of these error measures on the mesh nodal positions is accurately captured by adjoint sensitivity analysis. Additionally, a convolution-based gradient regularization strategy is used to ensure stable and effective adaptation of high-order meshes. We demonstrate that the proposed framework can improve mesh quality and reduce the error by up to 10 times for the solution of Poisson and linear elasto-static problems. The approach is general with respect to the dimensionality, the order of the mesh, the types of mesh elements, and can be applied to any PDE that admits well-defined adjoint operators.}

\keywords{Mesh optimization, $r$-adaptivity, high-order meshes, PDE-constrained optimization}
\maketitle

\section{Introduction}
\label{sec_intro}

Efficient and accurate solution of finite element (FE) discretizations is crucial for a wide range of scientific and engineering applications, including fluid dynamics and structural analysis. Specifically, high-order finite element (FE) methods have gained increasing popularity in recent years due to their scalability on modern architectures \cite{Fischer2002,fischer2020scalability,MFEM2020,CEED2021} and their superior accuracy compared to low-order methods. They are essential for achieving optimal convergence rates on curvilinear domains, preserving symmetry, and capturing key flow features in simulations with moving meshes \cite{Shephard2011, Dobrev2012, Boscheri2016}. The accuracy of an FE discretization directly depends on the quality and resolution of the mesh used to construct the finite element approximation space. Thus, deploying mesh adaptivity increases the accuracy of the FE solution while minimizing the increase in computational cost. Although traditional FE mesh adaptivity employs $h-$, $p-$, or $r$-adaptivity, we focus here on $r$-adaptivity for its ability to preserve element count, polynomial order, and mesh connectivity - making it especially well suited for dynamic simulations where mesh quality must adapt without altering the underlying discretization \cite{TMOP2020CAF}.

While it is not possible to determine the exact discretization error before solving the discretized PDE \cite{VERFURTH199}, various heuristic methods have been developed to estimate the error and adapt the mesh to improve solution accuracy. These approaches include \textit{a posteriori} error estimator methods \cite{baines1994moving,tang2005moving,Weizhang2011} and geometry or solution feature-based methods \cite{TMOP2021,TMOP2020CAF}. The general idea in these approaches is to derive a discrete monitor function that estimates the solution error, which is minimized throughout the mesh using $h$/$p$/$r$-adaptivity. Similarly, energy-based methods have been developed for non-linear PDEs by minimizing an objective function based on the potential energy of the mechanical system by $r$-adaptivity. The optimization reduces the discretization error measured by the energy norm \cite{Thoutireddy2004,Scherer2007}. Another category of methods is based on deriving a metric tensor field using goal-oriented error estimates that drive the anisotropic mesh adaptation  \cite{loseille2010fully,park2021verification,wallwork2022goal}.
Although all these approaches have been successfully applied to the $r$-adaptivity in various settings, they do not directly optimize the target quantity (monitor function, potential energy, error estimate) on the mesh resulting from optimization. Instead, they rely on the heuristic assumption that the adapted mesh will yield a good PDE solution, without explicitly controlling the target quantity on the new mesh. In this work, we aim to provide precisely that control.

We propose a novel PDE-constrained high-order mesh $r$-adaptivity framework that simultaneously improves mesh quality and solution accuracy. We formulate $r$-adaptivity as an optimization problem, where the objective is a convex combination of mesh quality metric and a measure based on solution error. We utilize the Target-Matrix Optimization Paradigm (TMOP) for the mesh quality component due to its ability to precisely prescribe geometric targets, its purely algebraic formulation applicable to various element types in both 2D and 3D \cite{TMOP2019SISC,TMOP2020CAF,TMOP2021EWC}, and its suitability for GPU acceleration \cite{camier2023accelerating}. The solution error measure is PDE-dependent and can incorporate existing a posteriori error estimators from the literature \cite{baines1994moving,tang2005moving,Weizhang2011}. In the context of elliptic PDEs, we explore error measures based on different norms of the discrete solution. The PDE-constrained $r$-adaptivity problem is solved using  techniques inspired by shape optimization \cite{HSU199,Swartz2023}. We solve the resulting optimization problem iteratively, using consistent gradients of the objective function with respect to mesh nodal positions, computed via adjoint sensitivity analysis. To mitigate localized gradients and numerical instabilities, we apply a convolution-based regularization filter \cite{Lazarov2011,Swartz2023}. A gradient-based optimizer then determines optimized nodal positions that minimize the objective of interest. As the proposed framework is intended for use in differentiable physics simulations and automated design optimization, it is designed to be general with respect to element types, spatial dimensions, mesh order, and any PDEs with well-defined adjoint operators.

The remainder of this paper is organized as follows: Section \ref{sec_prelims} introduces the PDEs of interest, our mesh representation, and the finite element discretization. Section \ref{sec::PDEconstOpt} provides an overview of PDE-constrained optimization in the context of $r$-adaptivity.  Sections~\ref{sec::PDEperformanceMeasure} and \ref{sec::qualityperformanceMeasure} describe the error measures and mesh quality metric, respectively, used to define the global objective. Section~\ref{sec::solver} discusses details of the implementation of the optimization solver.  Section~\ref{sec::AD} presents the use of automatic differentiation to calculate the gradients of the mesh quality metric. Numerical results are shown in Section~\ref{sec::numericalExamples}, and conclusions are summarized in Section~\ref{sec::Conclusion}.

\section{Preliminaries}
\label{sec_prelims}

The $r$-adaptivity approach introduced in this paper applies to any PDE with a well-defined adjoint operator. Without loss of generality, we demonstrate it on two representative elliptic problems: the Poisson equation and the linear elasto-static equation. Let $\Omega \subset \mathbb{R}^{d}$ be the physical domain with volume force $b$. The boundary of $\Omega$ consists of two parts, $\partial \Omega = \Gamma_D \cup \Gamma_N$, with Dirichlet boundary conditions enforced on $\Gamma_D$, and Neumann boundary conditions applied on $\Gamma_N$.

The governing equation for the Poisson problem is
\begin{subequations}
\label{eq::poisson_strong_form}
    \begin{align}
    \nabla ^2 u + b &= 0 \quad \text{in} \ \Omega,  \\
    u &= u_b \quad  \text{on} \ \Gamma_D, \\
    \nabla u \cdot \boldsymbol{n} &= t \quad  \text{on} \ \Gamma_N,
    \end{align}
\end{subequations}
where inhomogeneous boundary conditions are imposed on $\Gamma_D$ and $\Gamma_N$.

The governing equation for the considered linear elasto-static problem reads
\begin{subequations}
\label{eq::elasticity_strong_form}
    \begin{align}
    \nabla \cdot \bsig + \bb &= 0 \quad \text{in} \ \Omega,  \\
    \bu &= 0 \quad  \text{on} \ \Gamma_D, \\
    \bsig \cdot \boldsymbol{n} &= \mathbf{t} \quad  \text{on} \ \Gamma_N,
    \end{align}
\end{subequations}
where the Cauchy stress tensor is denoted as $\bsig(\bu) = \mathbf{D} \boldsymbol{\varepsilon}$, with $\mathbf{D}$ being the constitutive tensor, and $\boldsymbol{\varepsilon}$ the elastic strain tensor. The elastic infinitesimal strain $\boldsymbol{\varepsilon}$ is given by $\boldsymbol{\varepsilon} = \frac{1}{2} \left( \nabla (\boldsymbol{\mathit{u}}) + \nabla (\boldsymbol{\mathit{u}})^T \right)$.

The equations are discretized using finite elements on a mesh $\mcm$ resolving the physical domain $\Omega$. The accuracy of the finite element analysis of a boundary-value problem directly depends on the mesh quality and resolution used to construct the finite element approximation space. The $r$-adaptivity performed in this work aims to find node locations for a given mesh, optimizing the accuracy of a given discretized boundary value problem.

\subsection{Mesh representation}
In a finite element framework \cite{MFEM2020,MFEM2024}, the physical domain $\Omega \subset \mathbb{R}^{d}$, is discretized as a union of $N_E$ curved mesh elements, $\mcm^e$, $e=1\dots N_E$, each of order $p$. A discrete representation of these elements is obtained by selecting a set of scalar basis functions $\{ \phi_i \}_{i=1}^{N_p}$ on a reference element $\bar{\mcm}^e$. The bases, $\phi_i(\br)$, $\br \in \bar{\mcm}^e$, are usually chosen as Lagrange interpolation polynomials defined at a set of $N_p$ interpolation nodes $\zeta_j \in \br$, such that $\phi_i(\mathbf{\zeta}_j) = \delta_{ij}$. The nodal coordinates of an element $\mcm^e$ in the mesh $\mcm$ are described by a matrix $\bx^e$, of size $d \times N_p$. Note that for quadrilaterals/hexahedra, the reference element is $\bar{\mcm}^e \in [0,1]^d$ and $N_p = (p+1)^d$. For triangles/tetrahedra, the reference element is typically the standard simplex, with $N_p = \binom{p+d}{d}$.

Given $\bx^e$ and the choice of Lagrange interpolants as  bases, we introduce the map $\Phi_e:\bar{\mcm}^e \to \mathbb{R}^d$ whose image is the geometry of the physical element $\mcm^e$:
\begin{equation}
\label{eq_x}
\bx(\br) \vert_{\mcm^e} =
   \Phi_e(\br) \equiv
   \sum_{i=1}^{N_p} \bx^e_i \phi_i(\br),
   \qquad \br \in \bar{\mcm}^e, ~~ \bx=\bx(\br) \in \mcm^e,
\end{equation}
where $\bx^e_i$ denotes the coordinates of the $i$-the node of $\mcm^e$. Hereafter, $\bx(\br)$ will denote the element-wise position function defined by \eqref{eq_x} with $\mcm^e$ omitted for brevity.

Further, the Jacobian matrix $A$ represents the local deformation of $\mcm^e$ with respect to the $\bar{\mcm}^e$ at $\bbx$. The Jacobian of the mapping $\Phi_e$ at any reference point $\bbx \in \bar{\mcm}^e$ is
\begin{equation}
\label{eq_A}
  A_{ab}(\br) = \frac{\partial x_{a}(\br)}{\partial {r}_b} =
    \sum_{i=1}^{N_p} x^e_{i,a} \frac{\partial \phi_i(\br)}{\partial r_b}, \quad a,b = 1 \dots d
\end{equation}
where $x_a$ represents the $a$th component of $\bx$ \eqref{eq_x}, and $x^e_{i,a}$ represents the $a$th component of $\bx^e_i$, i.e., the $i$th DOF in element $\mcm^e$. This Jacobian enables the efficient use of the reference element for gradient calculation and integration. It can also be used to determine mesh validity, i.e., $\det(A)$ must be greater than 0 at every point in the mesh.

\subsection{Weak forms of the model problems}
\label{sec_weak}

For the Poisson problem the weak form is given by
\begin{equation}
\label{eq_poisson_weak0}
{\text{find }} {u}\in {\mcu}\quad {\text{such that }}
a\left({u},{v}\right)=l\left({v}\right)\quad \forall {v}\in\ {\mcu},
\end{equation}
where the solution and test spaces are $\mcu = H^{1}(\Omega)$, and the
bilinear $a\left(\cdot,\cdot\right)$ and linear $l\left(\cdot\right)$ forms are defined as
\begin{equation}
\label{eq_poisson_weak}
\begin{aligned}
a(u, v) &= \int_{\Omega} \nabla u : \nabla v \, d\Omega
        + \int_{\Gamma_D} \gamma\, u \cdot v \, d\Gamma, \\
l(v) &= \int_\Omega b \cdot v \, d\Omega
        + \int_{\Gamma_D} \gamma\, u_b \cdot v \, d\Gamma
        + \int_{\Gamma_N} t \cdot v \, d\Gamma.
\end{aligned}
\end{equation}
Inhomogeneous Dirichlet conditions are enforced weakly using the second term in the bilinear and linear forms. The penalty parameter $\gamma$ controls this weak enforcement and is typically set to a large constant, e.g., $10^5$ by default in the proposed implementation. The FE space used to discretize \eqref{eq_poisson_weak} is defined as
\begin{equation}
{\mcu}_h = \{u_h \in C^0(\Omega) ~ | ~
           u_h|_{\mcm^e} = \bar{u} \circ \Phi_e^{-1},
           \bar{u} \in \mathbb{Q}_{p_u}(\bar{\mcm}^e) \text{ or } \mathbb{P}_{p_u}(\bar{\mcm}^e) \},
\end{equation}
where the FE solution's polynomial order $p_u$ on each element can generally be different than the mesh order $p$; we use $p_u = p$ unless specified otherwise.

\noindent
The weak form for the considered linear elasticity problems is
\begin{equation}
{\text{find }} \boldsymbol{u}\in \boldsymbol{\mcu} \quad {\text{such that }}
a(\boldsymbol{u},\boldsymbol{v})=l(\boldsymbol{v}) \quad
\forall \boldsymbol{v} \in \boldsymbol{\mcu},
\end{equation}
where the solution and test spaces are defined as
\begin{equation}
\boldsymbol{\mcu} =
  \{\boldsymbol{u} \in \left[ H^{1}(\Omega) \right]^d,
  \boldsymbol{u} = 0 \ \mathrm{on} \ \Gamma_{D} \},
\end{equation}
and the bilinear $a\left(\cdot,\cdot\right)$ and linear $l\left(\cdot\right)$ forms are defined as
\begin{equation}
\label{eq_elastic_weak}
a(\boldsymbol{u},\boldsymbol{v}) =
  \int_{\Omega} \sigma(\boldsymbol{u}) : \nabla \boldsymbol{v} \ d\Omega, \quad
  l(\boldsymbol{v})=\int_\Omega b \cdot \boldsymbol{v} ~ d\Omega
	+ \int_{\Gamma_N} {t} \cdot \boldsymbol{v} ~ d\Gamma.
\end{equation}
The corresponding FE space to discretize \eqref{eq_elastic_weak} is
$\boldsymbol{\mcu}_h = \left[ ~ {\mcu}_h ~ \right]^d$.

\section[PDE-constrained r-adaptivity]{PDE-constrained $r$-adaptivity}
\label{sec::PDEconstOpt}

The proposed $r$-adaptivity approach aims to reduce a discretization error measure of the solution of a PDE by optimizing the placement of mesh nodes without changing the mesh connectivity and ensuring good element quality. To this end, we deploy gradient-based PDE-constrained optimization. Here, the optimization aims to minimize an objective function $\mcf$ based on a weighted combination of the solution error measure ($\mcf_p\left(\cdot,\cdot\right)$) and mesh quality ($\mcf_{\mu}\left(\cdot\right)$), subject to a constraint defined by a PDE. The problem is defined as
\begin{subequations}
\label{eq::linelast}
\begin{align}
    \underset{\boldsymbol{x}}{\min} ~~ &\mcf(\boldsymbol{x}) = \alpha \mcf_p\left(u(\boldsymbol{x}),\boldsymbol{x}\right) + \mcf_{\mu}\left(\boldsymbol{x}\right),    \\
     \text{subject to} \quad &\mcr_P\left(u; \boldsymbol{x}\right) = \boldsymbol{0}
\end{align}
\end{subequations}
where the PDE is written in residual form as $\mcr_{P}\left(u; \boldsymbol{x}\right)=a\left(u,v\right)-l\left(v\right) = 0$ and couples the state PDE solution $u$ and mesh control parameters, i.e., nodal coordinates $\boldsymbol{x}$. The coupling makes the optimization computationally challenging for large 3D problems. For simplicity, the formulation in \eqref{eq::linelast} uses the scalar PDE solution $u$ from \eqref{eq::poisson_strong_form}, though the method generalizes directly to vector-valued solutions $\boldsymbol{u}$, such as those arising in linear elasticity problems \eqref{eq::elasticity_strong_form}.

\begin{figure}[tbh!]
\begin{center}
$\begin{array}{cc}
\includegraphics[height=0.4\textwidth]{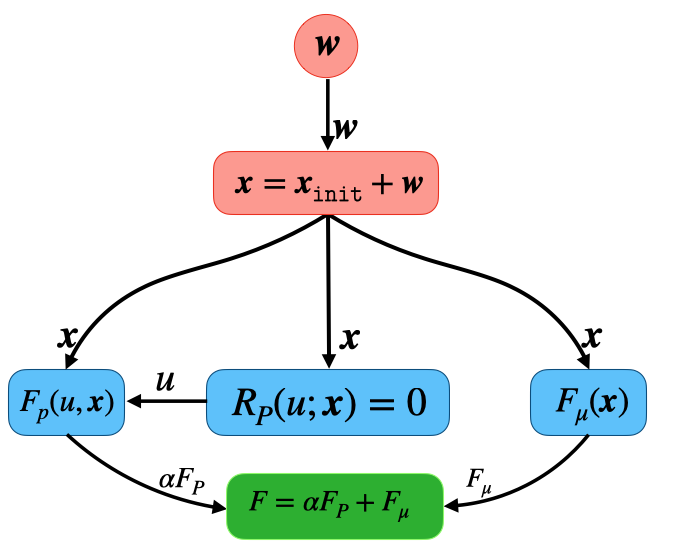} &
\includegraphics[height=0.4\textwidth]{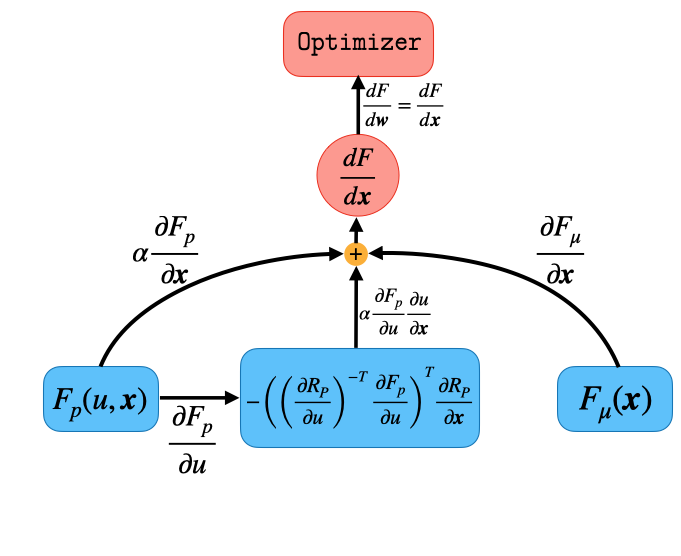} \\
\textrm{(a)} & \textrm{(b)} \\
\end{array}$
\end{center}
\caption{(a) \emph{Forward analysis} for computing the objective function and (b) \emph{reverse mode} for computing the gradient of the objective using the adjoint method.}
\label{fig_forwardgraph}
\end{figure}

We employ adjoint sensitivity analysis and gradient-based optimization to solve \eqref{eq::linelast} and determine optimized nodal locations on a given mesh, guided by the user-specified choice of $\mcf_p$, $\mcf_{\mu}$, and the PDE of interest. Figure \ref{fig_forwardgraph} gives a high-level overview of how the objective and its gradient are computed in the proposed framework.
In the so-called \emph{forward mode}, shown in Figure \ref{fig_forwardgraph}(a), the nodal positions $\bx$ are updated using the initial mesh ($\bxi$) and prescribed nodal displacements ($\boldsymbol{w}$) as $\boldsymbol{x} = \bxi + \bw$. The associated PDE is then solved ($\mcr_P=0$) at the updated positions, and the state solution $u$ is utilized to evaluate $\mcf_p$. In parallel, the updated mesh is used to assess the mesh quality using $\mcf_{\mu}$. Finally, the proposed objective is evaluated as $F=\alpha F_P + F_{\mu}$.

The mesh nodal displacements $\boldsymbol{w}$ are computed iteratively using gradient-based optimization, which requires the gradient of the objective function with respect to the control variables $\bw$.
Assuming $\mcf$ is differentiable and the PDE solution is continuously dependent on the mesh nodal displacements, the gradient of the objective is
\begin{eqnarray}
    \frac{d F}{d \bw} = \frac{d F}{d \bx} &=& \frac{\partial \mcf}{\partial \bx}+ \frac{\partial \mcf}{\partial u}\frac{\partial u}{\partial \bx},\\
    \label{eq_total_grad_obj2}
    &=& \frac{\partial \mcf_{\mu}}{\partial \bx}  +
        \alpha \frac{\partial \mcf_P}{\partial \bx} + \alpha \frac{\partial \mcf_P}{\partial u}\frac{\partial u}{\partial \bx} .
\end{eqnarray}
The first two terms on the right-hand side of \eqref{eq_total_grad_obj2} capture the explicit dependence of the objective on mesh nodal positions, while the third term captures the implicit dependence through the PDE solution $u$. This gradient computation is illustrated in Figure \ref{fig_forwardgraph}(b) and described in more detail next.

\subsection{Adjoint sensitivity analysis}
\label{sec::ShapeAdjoint}

The explicit dependence of the objective on the nodal positions is straightforward to compute using direct differentiation based on the choice of $F_P$ and $F_{\mu}$.
In contrast, the implicit dependence of the objective on mesh nodal positions requires $\partial u/\partial \bx$. This dependence is computed using adjoint sensitivity analysis, which is a standard technique in shape and topology optimization frameworks \cite{Allaire2004,Bendsoe2004}.

As the state variable field $u$ satisfies the residual equation $\mcr_{P}(u;\bx)=0$
\begin{eqnarray}
    \frac{d \mcr_P}{d \bx} = & \frac{\partial \mcr_P}{\partial \bx} +
    \frac{\partial \mcr_P}{\partial u}
    \frac{\partial u}{\partial \bx} = 0,\\
    \label{eq:residual_mat_der_2}
    \implies & \frac{\partial u}{\partial \bx} = - \bigg(\frac{\partial \mcr_P}{\partial u}\bigg)^{-1}  \frac{\partial \mcr_P}{\partial \bx},
\end{eqnarray}
where the Jacobian $\frac{\partial \mcr_{P}}{ \partial u}$ is the stiffness matrix (tangent) used to solve the associated PDE, and thus readily available. $\frac{\partial \mcr_P}{\partial \bx}$ depends on the PDE of choice and is obtained via direct differentiation.

Following \eqref{eq_total_grad_obj2} and \eqref{eq:residual_mat_der_2},
the sensitivity $\frac{\partial  u}{\partial \bx}$ is removed from Eq. (\ref{eq_total_grad_obj2}) by first evaluating the adjoint variable $\lambda_u$ that solves
\begin{equation}
\label{eq:adjoint_solve}
    \lambda_u = \left( \frac{\partial \mcr_{P}}{ \partial u} \right)^{-T}\frac{\partial \mcf_P}{\partial u},
\end{equation}
and then computing
\begin{equation}
\frac{\partial \mcf_P}{\partial u}\frac{\partial u}{\partial \bx} = \lambda_u^T \frac{\partial \mcr_{P}}{\partial \bx} \equiv -{\bigg( \left( \frac{\partial R_P}{\partial u} \right)^{-T} \frac{\partial F_p}{\partial u} \bigg)}^T \frac{\partial R_P}{\partial \boldsymbol{x}}.
\label{RedWRTDesigVariable}
\end{equation}
In \eqref{eq:adjoint_solve}, $\frac{\partial \mcf_P}{\partial u}$ is the \emph{adjoint load} that depends on $F_P$.
This gradient computation using adjoint sensitivity analysis is illustrated in Figure \ref{fig_forwardgraph}(b).

\subsection{Convolution filter}

\begin{figure}[t!]
\begin{center}
$\begin{array}{cc}
\includegraphics[height=0.4\textwidth]{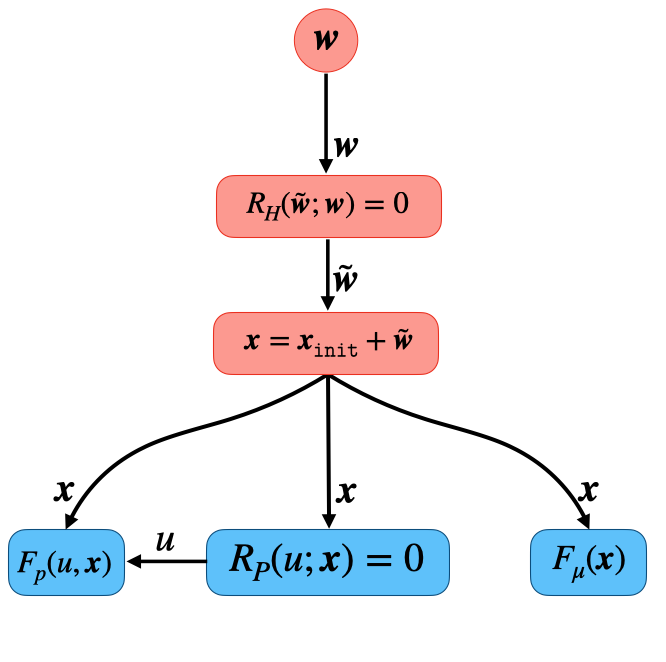} &
\includegraphics[height=0.4\textwidth]{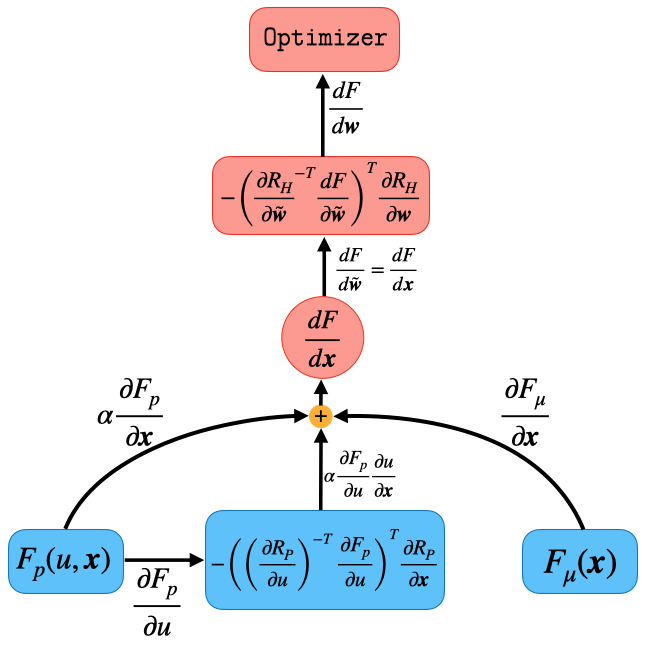} \\
\textrm{(a)} & \textrm{(b)} \\
\end{array}$
\end{center}
\caption{(a) \emph{Forward analysis} for computing the objective function and (b) \emph{reverse mode} for computing the gradient of the objective using the adjoint method, with convolution filter used to regularize the prescribed nodal displacements $\bw$.}
\label{fig_forwardgraphHelmholtz}
\end{figure}

It is well established that the usage of point-wise node displacements $\bw$ results in localized gradients, which can lead to numerical instabilities in the optimization process \cite{Lazarov2011}. In the context of $r$-adaptivity, localized gradients of the solution-based error measure can lead to severe mesh distortion. This issue is further exacerbated for high-order meshes. While the mesh quality term $\mcf_{\mu}$ in \eqref{eq::linelast} somewhat counters this issue, this term would become active only after the mesh is already distorted by the instabilities. We employ a convolution-type filter widely utilized in topology optimization frameworks \cite{Lazarov2011} to compute a continuous nodal displacement field $\tbw$ with a certain regularity, which is controlled by the filter radius $\delta$. The convolution filter defines the filtered node displacement field $\tbw$ as a solution of a vector PDE filter equation; that is, we compute
\[
\tbw \in \mathcal{H}_{\mathrm{H}} =
\left\{ \tbw \in \left[H^{1}\! \left( \Omega \right) \right]^d \right\}
\]
that satisfies
\begin{equation}
\label{Eq_Helmholtz}
a\left(\tbw,\bmiv\right) = l\left(\bw; \bmiv\right)
\end{equation}
for all $\bmiv \in \mathcal{H}_{\mathrm{H}}$, where
\begin{align}
   a\left(\tbw,\bmiv\right) &=  \int_{\Omega^d} \left(\ \delta^2 \nabla { \tbw} : \nabla  \bmiv + {\tbw} \cdot \bmiv \right) \ d \Omega \\
   l\left(\bw; \bmiv\right) &= \int_{\Omega^d}  {\bw}  \cdot \bmiv  \ d \Omega.
\end{align}
The parameter $\delta$ controls the smoothness of the filtered nodal displacement field $\tbw$. These filtered displacements are then used to compute the mesh nodal positions as $\bx = \bxi + \tbw$. This modification is shown in Figure \ref{fig_forwardgraphHelmholtz}, where the original forward-mode graph is extended by inserting the convolution filter $\mcr_{H}$, where $\mcr_{H}=a\left(\tbw,\bmiv\right)-l\left(\bw; \bmiv\right)$. Note that this entails an additional adjoint solution for the computation of the non-filtered sensitivities, as shown in Figure \ref{fig_forwardgraphHelmholtz}(b).

\section{Error measures for the PDE solution}
\label{sec::PDEperformanceMeasure}

While the method supports any differentiable error measure, we select three representative choices for $\mcf_p$ to demonstrate the optimization procedure. The first relies on a traditional error indicator commonly used in mesh refinement, while the others, though not classical error measures, are designed to indirectly reduce the discrepancy between the discrete and exact PDE solutions for smooth problems. Note, each of the error measures presents a particular choice of mesh optimization. The optimal solution for one of the measures is thus not necessarily the optimal solution for any of the other ones.

\subsection{Element-local variation}
\label{sec::perf_sol_local}
The element-local variation measure compares the local solution with the element-wise mean:
\begin{equation}
\mcf_P(\bx,{u}(\bx)) = \int_{\Omega} ({u}- \mathcal{B}({u}))^2 \ d  \Omega,
\label{eq::errorSolMean}
\end{equation}
with scalar $\mathcal{B}({u})$ defined as
\begin{equation}
\mathcal{B}({u})|_{\Omega^e} = \frac{\int_{\Omega^e}u}{|\Omega^e|}.
\label{eq::sol_mean}
\end{equation}
With \eqref{eq::errorSolMean}, the mesh becomes denser in areas with steep physical gradients and sparser in areas with gradients closer to zero. This error measure is inspired by the work of Zahr et al. \cite{zahr2018optimization}.

\subsection{Indirect minimization of the energy norm of the error}
\label{sec::compliance}

Let ${u}\in{\mcu}$ denote the exact continuous and ${u}_h\in{\mcu}_h$ the exact discrete solution of the associated PDE, where ${\mcu}_h \subset {\mcu}$ is a standard finite element space defined on a given mesh $\mcm$ \cite{Braess2007}. Following \eqref{eq_poisson_weak0}-\eqref{eq_poisson_weak}, the energy norm of the error ${e}_h={u}-{u}_h$ is given as
\begin{align}
    a\left({e}_h,{e}_h\right)=a\left({u}-{u}_h,{u}-{u}_h\right).
\end{align}
Using the Galerkin orthogonality of the discrete solution with respect to the complementary space ${\mcu}/{\mcu}_h$, the above equality can be written as
\begin{align}
\label{eq_en100}
    a\left({u}-{u}_h,{u}-{u}_h\right)=a\left({u},{u}\right)-a\left({u}_h,{u}_h\right)=l\left({u}\right)-l\left({u}_h\right),
\end{align}
where ${u}_h$ solves exactly $a\left({u}_h,{v}_h\right)=l\left({v}_h\right),\ \forall {v}_h\in {\mcu}_h$. Since the exact energy $a\left({u},{u}\right)$ in \eqref{eq_en100} does not depend on the mesh $\mcm$, maximizing the \emph{load functional} $l\left({u}_h\right)$ is equivalent to minimizing the energy norm of the error $a\left({e}_h,{e}_h\right)$.

The proposed objective is then
\begin{equation}
\mcf_P(\bx,{u}(\bx)) =  -l\left({u_h}\right)=-\int_{\Omega} {u_h} \cdot {f} \ d  \Omega,
\label{eq::errorCompl}
\end{equation}
where the negative sign is introduced in front of the load functional such that minimizing $\mcf_P$ minimizes the energy norm of the error.
Note that \eqref{eq::errorCompl} is a common objective for topology and shape optimization frameworks \cite{Bendsoe2004, Allaire2004}, where the load functional with positive sign (also known as compliance) is minimized for a given constraint on the volume of material. This minimization process ensures that the final shape or the material distribution provides the stiffest possible design among all possible designs with a given amount of material.

An equivalent formulation is presented in \cite{Delfour1985}, where the problem is formulated by considering an energy functional associated with the elliptic problem. The minimum of this functional corresponds to the optimal solution; however, as discussed in \cite{Delfour1985}, the non-convexity with respect to the nodal positions and the fixed mesh topology restricts the optimization process to finding an adapted mesh that minimizes energy error without guaranteeing a global optimum.

\subsection{Gradient continuity}
\label{sec::perf_zz}
\begin{figure}[t!]
\begin{center}
$\begin{array}{cc}
\includegraphics[height=0.6\textwidth]{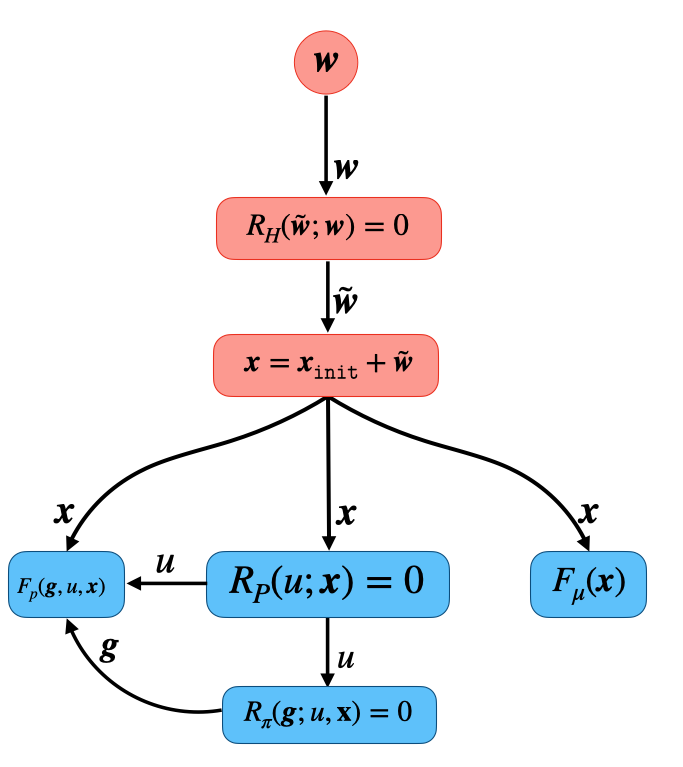} &
\includegraphics[height=0.6\textwidth]{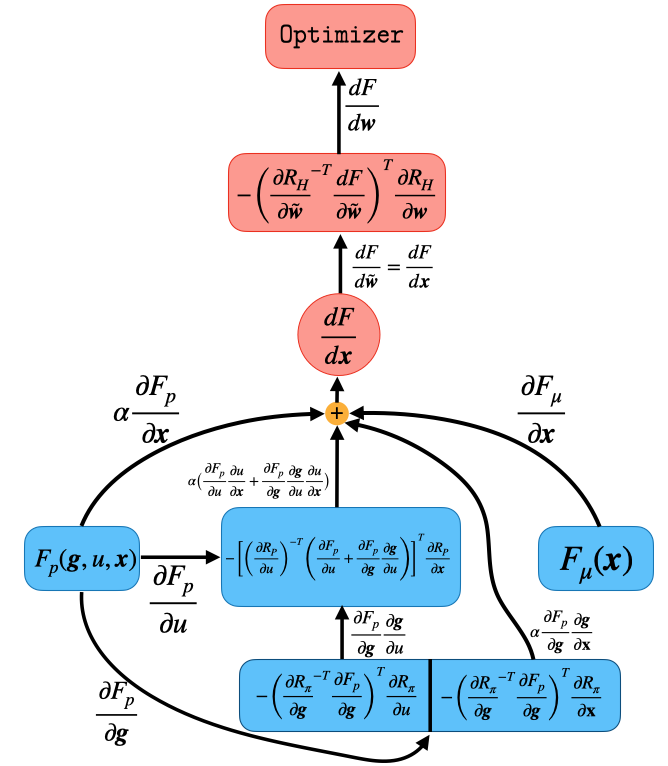} \\
\textrm{(a)} & \textrm{(b)} \\
\end{array}$
\end{center}
\caption{(a) \emph{Forward analysis} for computing the objective function and (b) \emph{reverse mode} for computing the gradient of the objective using the adjoint method, for the gradient continuity measure.}
\label{fig_forwardgraphZZ}
\end{figure}

The third error measure is inspired by the Zienkiewich-Zhu error estimator \cite{zienkiewicz1987simple}. It uses the mismatch in the gradient of the solution across element boundaries as a measure for the error in the solution. Specifically, it measures the $L^2$ norm between the gradient of the solution ${\nabla u}$ and a globally smoothed gradient field ${g}$ as:
\begin{equation}
\mcf_P(\bx,u({\bx})) =
\int_{\Omega} ({\nabla u}-  {g})^2 \ d  \Omega,
\label{eq::errorGrad}
\end{equation}
where the globally smoothed gradient field
${g} \in \left[H^{1} \left( \Omega \right)\right]^d$ is computed by the following projection:
\begin{equation}
\mcr_{{\pi}}({g};\bx,u(\bx)) =
\int_{\Omega}  ({g}-\nabla u)\cdot v \ d  \Omega = 0
\label{eq::GlobalZZ}
\end{equation}
for all $v \in \left[H^{1} \left( \Omega \right)\right]^d$. The sensitivity analysis must now account for additional implicit dependencies introduced in $\mcf_P$:
\begin{equation}
    \frac{d F_p}{d \boldsymbol{x}} = \frac{\partial F_p}{\partial \boldsymbol{x}} + \frac{\partial F_p}{\partial u}\frac{\partial u}{\partial \boldsymbol{x}}+
    \frac{\partial F_p}{\partial {g}}\bigg(\frac{\partial {g}}{\partial \boldsymbol{x}} + \frac{\partial {g}}{\partial u}\frac{\partial u}{\partial \boldsymbol{x}}\bigg),
\end{equation}
similar to the procedure in Section \ref{sec::ShapeAdjoint}, as illustrated in Figure \ref{fig_forwardgraphZZ}.
With this error measure, the mesh is expected to become denser in areas with large changes in the second derivative of the physical field and sparser in areas with constant gradients.

\section{Mesh quality measure}
\label{sec::qualityperformanceMeasure}

The mesh quality measure computation is based on the TMOP approach for high-order meshes \cite{TMOP2019SISC}. TMOP enables control over local mesh quality, computed through the local Jacobian matrix $A_{d \times d}$ defined in \eqref{eq_A}, while still optimizing the mesh globally.

The first step for mesh quality improvement is the specification of a target matrix $W_{d \times d}$, analogous to $A_{d \times d}$, that prescribes the target deformation at each point in the mesh. Target matrix construction is typically driven by the desired geometric properties of the mesh, as any Jacobian matrix can be written as a composition of the four geometric components, namely volume, rotation, skewness, and aspect-ratio:
\begin{equation}
\label{eq_W}
W_{d\times d} = \underbrace{\zeta}_{\text{[volume]}} \circ
    \underbrace{R_{d\times d}}_{\text{[rotation]}} \circ
    \underbrace{Q_{d\times d}}_{\text{[skewness]}} \circ
    \underbrace{D_{d\times d}}_{\text{[aspect-ratio]}}.
\end{equation}
    In the current work, we are interested in ensuring good element \emph{shape} (skewness and aspect-ratio), so we prescribe $W$ to be that of an ideal element, i.e., square for quad elements, cube for hex elements, and equilateral simplex for triangles and tetrahedra. The reader is referred to \cite{knupp2019target,TMOP2019SISC,TMOP2020CAF,TMOP2021EWC} for advanced techniques on target construction for geometry- and simulation-driven adaptivity.

\begin{figure}[tb]
\centerline{\includegraphics[width=0.4\textwidth]{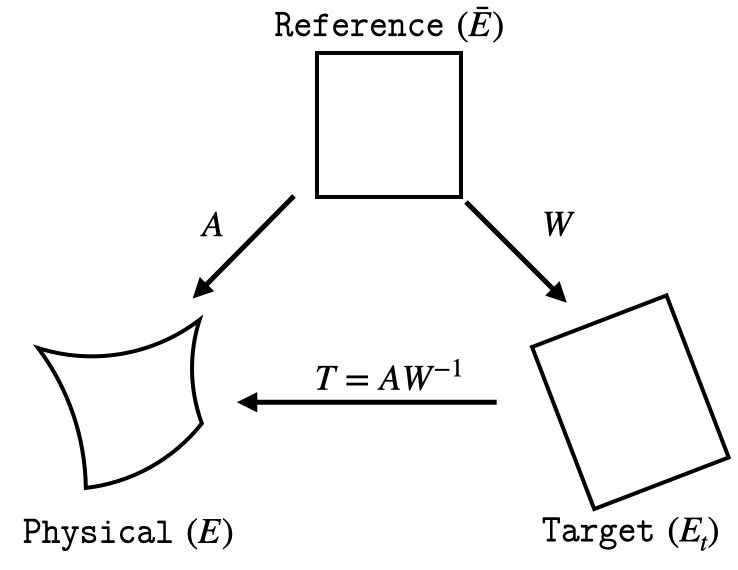}}
\vspace{-1mm}
\caption{Schematic representation of the major TMOP matrices.}
\label{fig_tmop}
\vspace{-3mm}
\end{figure}

Next, a mesh quality metric is chosen to measure the deviation between the current ($A$) and target Jacobians ($W$) in terms of one or more geometric parameters. These mesh quality metrics are typically defined using the target-to-physical transformation $T=AW^{-1}$ shown in Figure \ref{fig_tmop}. For example, $\mu_{2,s}(T)=\frac{\mid\mid T \mid\mid^2}{2\tau}-1$ is a \emph{shape} metric\footnote{The metric subscript follows the numbering in \cite{Knupp2020, knupp2022geometric}.} that depends on the skewness and aspect-ratio components, but is invariant to orientation/rotation and volume. Here, $\mid\mid T\mid\mid$ and $\tau$ are the Frobenius norm and determinant of $T$, respectively.
Seldom, we also use metrics that directly depend on $A$ and $W$ (instead of via $T$) and we denote them as $\nu(A,W)$. For example, $\nu_{107,OS}(A,W)= 0.5 ({\tt det}\,A)^{-1} \mid\mid A - (\mid\mid A \mid\mid/\mid\mid W\mid\mid) W \mid\mid^2$ is a \emph{shape+orientation} metric that depends on the skewness, aspect-ratio, and orientation components but is invariant to the volume of the element.

The quality metric $\mu(T)$ (or $\nu(A,W)$) leads to the following measure used for \ra
\begin{equation}
\label{eq_F_full}
  F_{\mu}(\bx) = \sum_{e = 1}^{N_E} F_{\mu,\mcm^e}(\bx^e)
               = \sum_{e = 1}^{N_E} \int_{\mcm^e_t} \mu(T(\bx^e)),
\end{equation}
where $F_{\mu}$ is the sum of element-wise quality
and $\mcm^e_t$ is the target element corresponding to the element $\mcm^e$. In prior work, we have used \eqref{eq_F_full} and its extension to demonstrate geometry and simulation-driven $r$-adaptivity in two- and three-dimensions for different element types \cite{TMOP2019SISC,TMOP2020CAF,TMOP2021EWC,barrera2023high}. Here, we aim to optimize for a general error measure based on the PDE solution with $F_{\mu}$ ensuring desired element quality.

\section{Optimization solver}
\label{sec::solver}

Optimal node locations can be determined by solving the proposed PDE-constrained optimization problem \eqref{eq::linelast} using standard techniques such as adjoint-based first-order gradient methods, sequential quadratic programming, augmented Lagrangian methods, or derivative-free approaches. We deploy a first-order gradient-based method with a ``black box'' optimizer as it provides a good balance between speed, accuracy, and implementation complexity. More specifically, to minimize the global objectives in \eqref{eq::linelast}, we use the Method of Moving Asymptotes (MMA) \cite{Svanberg2002}, a widely adopted first-order optimization algorithm in topology optimization frameworks. Although other first-order methods, such as gradient descent and the Broyden–Fletcher–Goldfarb–Shanno (BFGS) algorithm, could be used, we specifically choose MMA for its robustness and efficiency in solving large-scale problems with numerous design variables.

At each optimization iteration, the gradient of the objective is computed using the sensitivity analysis in Section \ref{sec::ShapeAdjoint} with gradients provided as input to the MMA solver. The MMA solver updates the unfiltered displacement vector $\bw$ that is then used in a standard line search to ensure that the resulting mesh ($\bx = \bxi + \tbw$) is valid. Mesh validity is checked using the determinant of the Jacobian of the transformation \eqref{eq_A} at a set of quadrature points in each element. If the minimum determinant of the Jacobian in the mesh is non-positive, we iteratively scale the displacement vector $\tbw$ by $0.5$ until the minimum determinant becomes positive. Optionally, we can include additional line-search constraints to ensure that the global objective and the norm of its gradient decrease sufficiently following each MMA iteration. In this work we use a maximum of 300 optimization iterations or until the ratio of the norm of the gradient of the objective with respect to the current and initial mesh decreases under a certain threshold. In future work, we will explore more convergence criteria that are better suited for the proposed formulation.

The linear system required for the forward and adjoint problems is solved using the standard Conjugate Gradient method with AMG-based preconditioning. The proposed methodology is verified using the high-order finite element library MFEM \cite{MFEM2024} with results presented in Section \ref{sec::numericalExamples}.

\section{Automatic differentiation for mesh quality computations}
\label{sec::AD}

The use of gradient-based optimization methods to minimize the objective requires computation of the derivatives of $\mcf_P$ and $\mcf_{\mu}$ with respect to the nodal positions. The explicit contribution from the mesh quality measure requires the derivative of the mesh quality metric (see Section 3.3 of \cite{camier2023accelerating}). Since these nonlinear metrics include matrix functionals,  they can lead to fairly complex expressions. As a result, prototyping new metrics and their derivatives can be time-consuming and error-prone. Thus, to simplify the implementation complexity and speed up the introduction of new mesh quality measures, we deploy automatic differentiation (AD) based on dual numbers and natively implemented in the MFEM discretization library \cite{MFEM2024,Andrej2025}.

Automatic differentiation (AD) with dual numbers is a powerful technique for computing exact derivatives numerically without relying on symbolic algebra or finite-difference based approximations. In this approach, the inputs $x \in \mathbb{R}$ are replaced with the corresponding dual numbers $\breve{x}=\langle x,x' \rangle$, where $\{x,x'\} \in \mathbb{R}$, and the functions $G(x)$ are replaced with \emph{dual-functions} (DF) $\mathcal{G}(\breve{x})=\langle G(x),x'G_x'(x)\rangle $, where $G_x'(x)$ denotes the derivative of $G(x)$ with respect to $x$.  Similarly, multivariate functions $G(\mathbf{x}):=G(x_1,x_2,\dots,x_N)$ are replaced with DFs $G(\mathbf{\breve{x}}) = \langle G(\mathbf{x}),\sum_i x_i'G_{x_i}'(\mathbf{x})\rangle $. AD is achieved by defining DFs corresponding to the primitive operators (e.g., addition, multiplication, trigonometric functions) that combine to form the target functional. For example, the DFs corresponding to subtraction and multiplication are
\begin{eqnarray}
\label{eq_ad_1}
    \breve{u}-\breve{v}=\langle u,u'\rangle  - \langle v,v'\rangle  = \langle u-v,u'-v'\rangle ,\\
\label{eq_ad_2}
    \breve{u}\cdot\breve{v}=\langle u,u'\rangle \cdot \langle v,v'\rangle  = \langle uv, uv'+u'v\rangle ,
\end{eqnarray}
and they combine for the determinant of a matrix $T=[\breve{a},\breve{b};\breve{c},\breve{d}]$ as
\begin{eqnarray}
\label{eq_ad_det}
{\tt det}(T) = \breve{a}\cdot\breve{d}-\breve{b}\cdot\breve{c}:=\underbrace{\langle ad-bc,ad'+da'-bc'-cb'\rangle }_{\text{following } \eqref{eq_ad_1}-\eqref{eq_ad_2}}.
\end{eqnarray}
Note that the last term in \eqref{eq_ad_det} is implicitly obtained through the subtraction and multiplication DFs in \eqref{eq_ad_1} and \eqref{eq_ad_2}, respectively. Finally, setting the dual component of an argument to 1 yields the derivative of the target function with respect to that argument in the dual component of the function's return value. For example, setting $a'=1$ and $b'=c'=d'=0$ in \eqref{eq_ad_det} returns the determinant value $ad-bc$ in the primal part and $d$ in the dual part, which matches the matrix determinant ${\tt det}(T)=ad-bc$, and its derivative w.r.t. $a$: $\partial {\tt det}(T)/\partial a=d$.

Similarly, other matrix functionals can be differentiated using AD. For example, AD of $\nu_{107}$ requires the primitive operations addition and multiplication for computing the Frobenius norm of $A$ and $W$, subtraction and multiplication for computing the determinant of $A$, and division to compute the ratio of the Frobenius norms of $A$ and $W$. AD can be deployed in the evaluation of gradients of the error measure, i.e. $\mcf_p$ in \eqref{eq::linelast} with preliminary results in adjoint-based gradient evaluation shown in \cite{Andrej2025}.

\section{Numerical results}
\label{sec::numericalExamples}
We present various numerical experiments to demonstrate the effectiveness of the proposed method and the impact of the selected error measure. All presented experiments were implemented using MFEM \cite{MFEM2024}.

\subsection{Mesh adaptation for 2D and 3D Poisson problems}

\begin{figure}[htbp]
\begin{center}
$\begin{array}{cccc}
\includegraphics[height=0.24\textwidth]{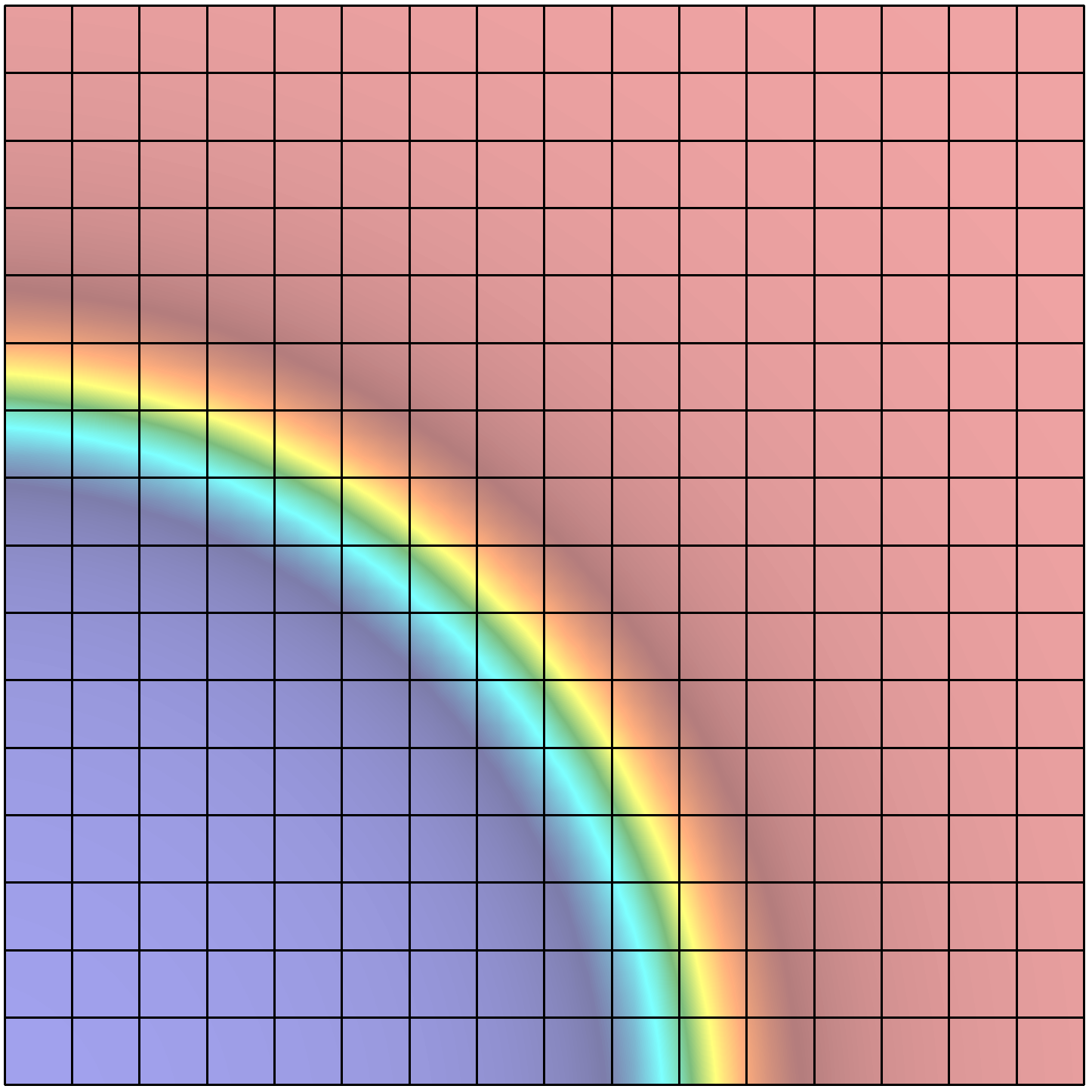} &
\includegraphics[height=0.24\textwidth]{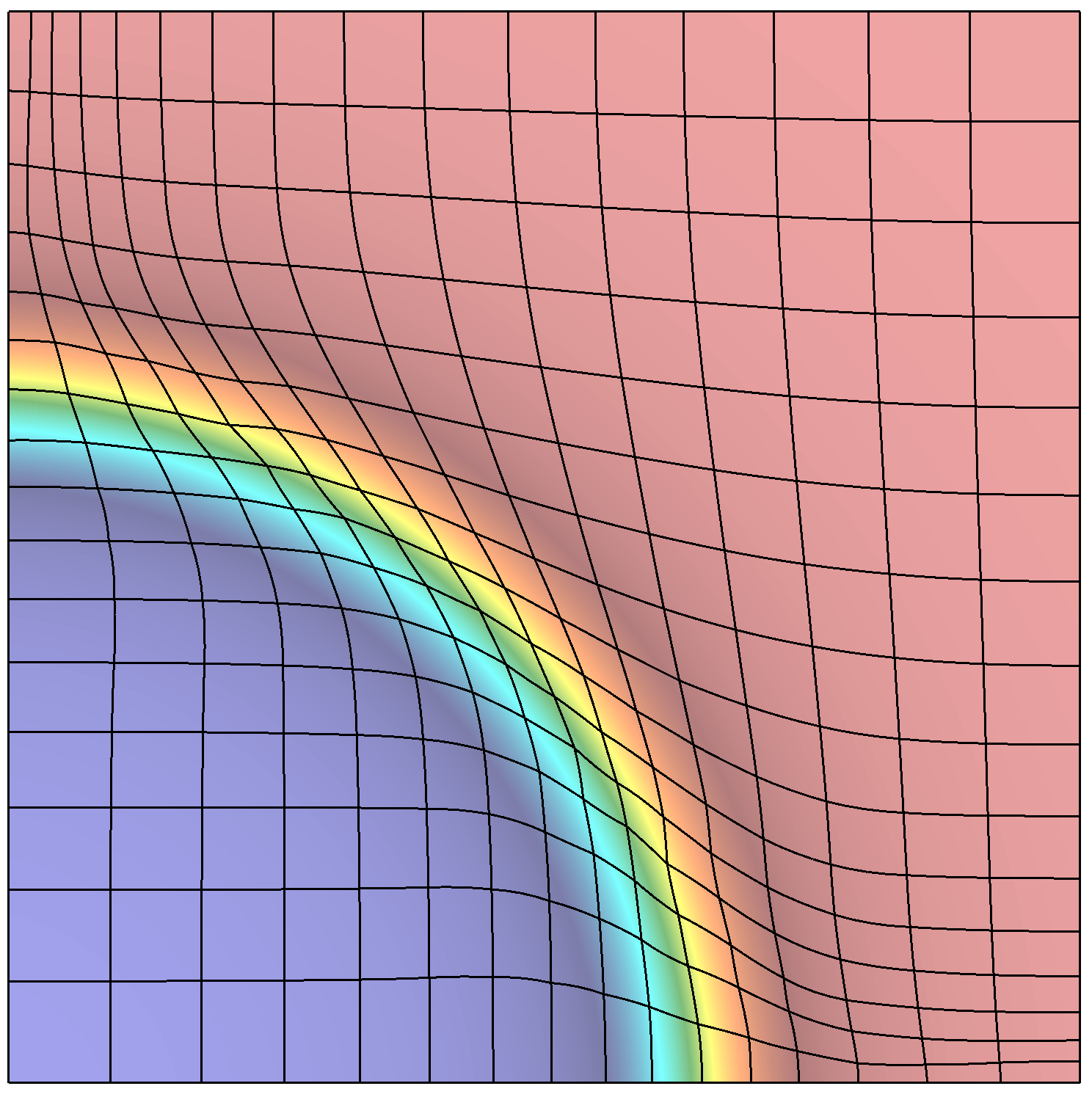} &
\includegraphics[height=0.24\textwidth]{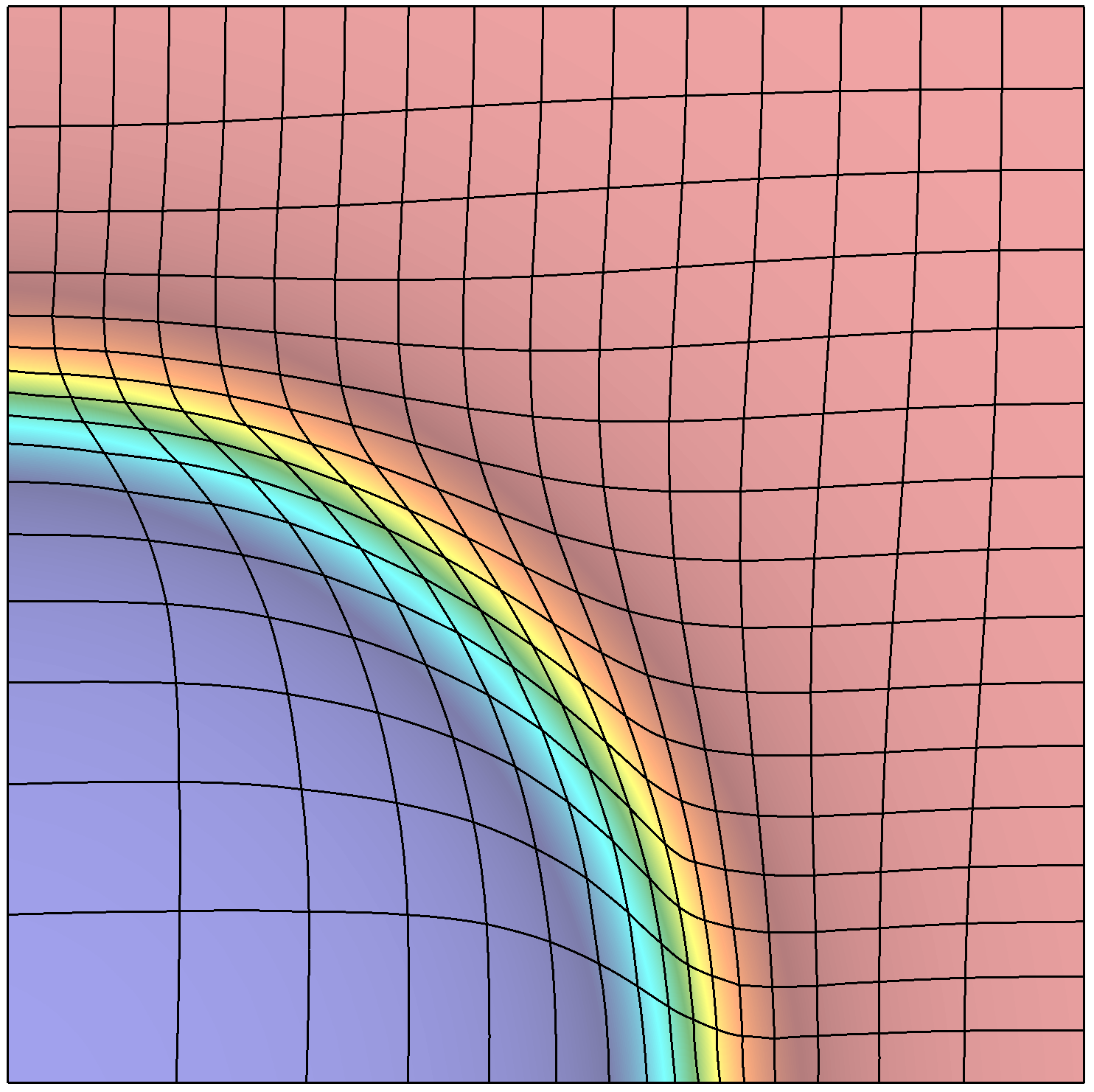} &
\includegraphics[height=0.24\textwidth]{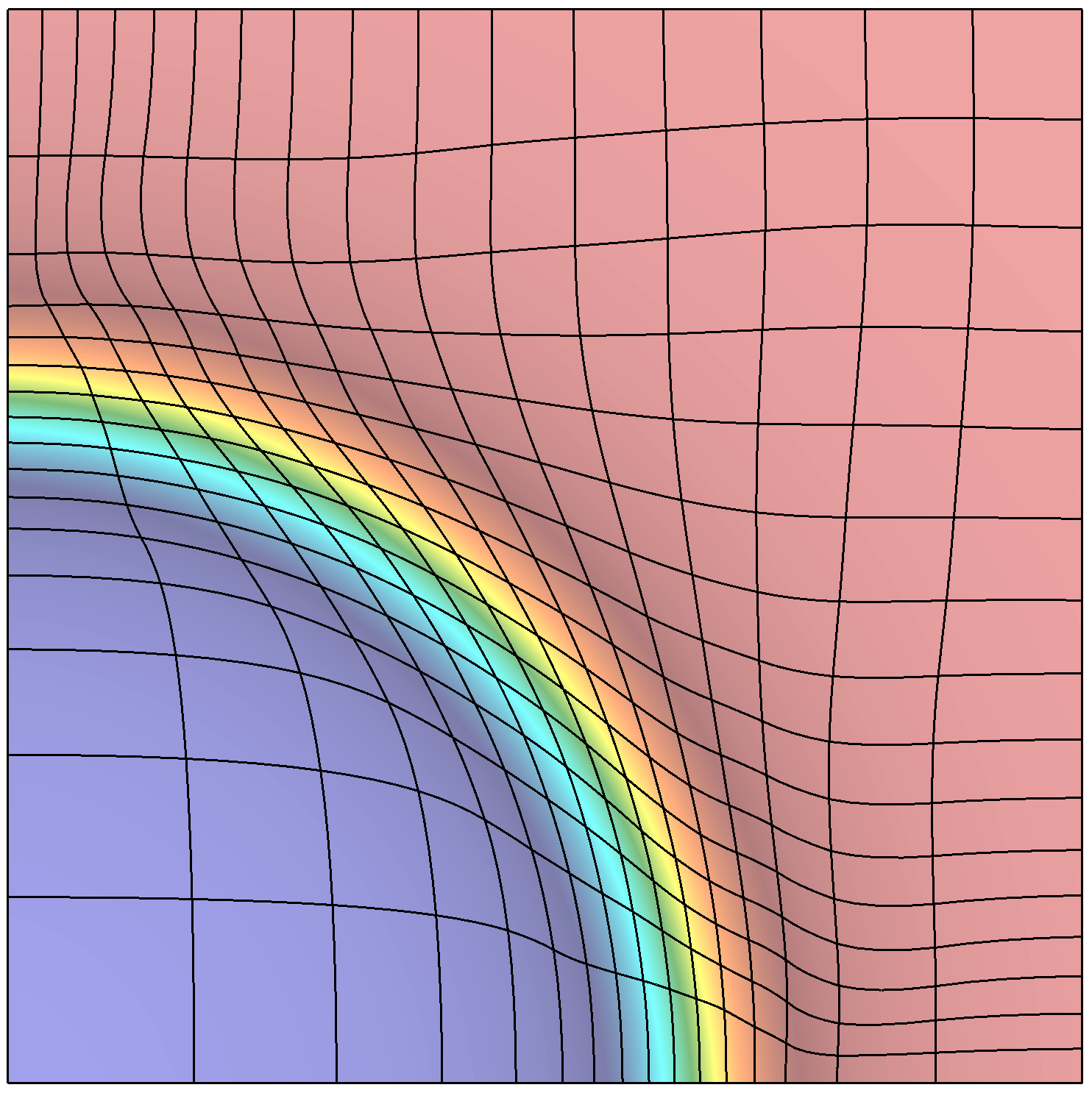} \\
\textrm{(a)} & \textrm{(b)} & \textrm{(c)} & \textrm{(d)}
\end{array}$
\end{center}
% \vspace{-7mm}
\caption{PDE-constrained optimization of a uniform $16\times16$ quadratic mesh. (a) Exact solution of the Poisson PDE \eqref{eq_exact_poiss_1} on the initial mesh, and meshes optimized with (b) the gradient continuity error measure (c) the element-local solution variation error measure and (d) the load functional-based error measure.}
\label{fig_poisson_circular}
\end{figure}

This first example demonstrates the \ra\ functionality of the presented method by approximating an analytical solution field $u(\bx)$ shown in Figure \ref{fig_poisson_circular}(a):
\begin{subequations}\label{eq_exact_poiss_1}
\begin{align}
    -\nabla^2 u &= f, \\
    u(\bx) &= \tan^{-1}\left(\beta\left(r(\bx) - r_c\right)\right),\\
    r(\bx) &=
\begin{cases}
\sqrt{(x - x_c)^2 + (y - y_c)^2} & \text{if } \mathbf{x}\in\mathbb{R}^2 \\
\sqrt{(x - x_c)^2 + (y - y_c)^2 + (z - z_c)^2} & \text{if } \mathbf{x}\in\mathbb{R}^3
\end{cases}
\end{align}
\end{subequations}
where $x_c = y_c = z_c = -0.05$, $r_c=0.7$, $\beta=20$. Inhomogeneous Dirichlet conditions are imposed on all boundaries.

We adapt linear, quadratic, and cubic two-dimensional quadrilateral meshes in $[0,1]^d$ to demonstrate that the proposed method recovers the expected convergence for the $L_2$ error in the solution. For this convergence study, a uniformly-sized Cartesian-aligned $4\times 4$ mesh is refined up to 3 times before the PDE-constrained optimization problem is solved with each of the three error measures in Section \ref{sec::PDEperformanceMeasure}. These error measures are well-suited for the solution of \eqref{eq_exact_poiss_1} as it has high localized gradients. For the mesh quality measure, the shape metric $\mu_2$ is used with the target matrix $W=I$ to encode the target element as an ideal square.

\begin{figure}[htbp]
  \centering
  \begin{subfigure}[b]{0.48\textwidth}
    \centering
    \resizebox{1.0\textwidth}{!}{%
    \begin{tikzpicture}
      \begin{axis}[
          legend columns=2,
          xlabel={$h$},
          ylabel={$||u_h-u^*||_{L_2}$},
          legend pos=north east,
          x dir=reverse,
          xmode=log,
          ymode=log,
          log basis x=10,
          log basis y=10,
          ymax=1e2,
          ymin=1e-8,
          legend cell align={left},
          every axis plot/.append style={line width=1pt}
        ]
        \addplot[color=red,mark=o,solid] table [x index=3, y index=5] {figures/tikz/avg_poisson_o1.txt};
        \addplot[color=red,dashed,mark=o,mark options={solid}] table [x index=3, y index=9] {figures/tikz/avg_poisson_o1.txt};
        \addplot[color=blue,mark=square,solid,mark options={solid}] table [x index=3, y index=5] {figures/tikz/avg_poisson_o2.txt};
        \addplot[color=blue,mark=square,dashed,mark options={solid}] table [x index=3, y index=9] {figures/tikz/avg_poisson_o2.txt};
        \addplot[color=black!60!green,mark=square,solid,mark options={solid}] table [x index=3, y index=5] {figures/tikz/avg_poisson_o3.txt};
        \addplot[color=black!60!green,mark=square,dashed,mark options={solid}] table [x index=3, y index=9]  {figures/tikz/avg_poisson_o3.txt};
        \addlegendentry{Initial ($p=1$)};
        \addlegendentry{Optimized ($p=1$)};
        \addlegendentry{Initial ($p=2$)};
        \addlegendentry{Optimized ($p=2$)};
        \addlegendentry{Initial ($p=3$)};
        \addlegendentry{Optimized ($p=3$)};
      \end{axis}
    \end{tikzpicture}
    }
    \caption{Element-local variation measure}
    \label{fig:avg_quad_l2}
  \end{subfigure}
  \hfill
  \begin{subfigure}[b]{0.48\textwidth}
    \centering
    \resizebox{1.0\textwidth}{!}{%
    \begin{tikzpicture}
      \begin{axis}[
          legend columns=2,
          xlabel={$h$},
          ylabel={$||u_h-u^*||_{L_2}$},
          legend pos=north east,
          x dir=reverse,
          xmode=log,
          ymode=log,
          log basis x=10,
          log basis y=10,
          ymax=1e2,
          ymin=1e-8,
          legend cell align={left},
          every axis plot/.append style={line width=1pt}
        ]
        \addplot[color=red,mark=o,solid] table [x index=3, y index=5] {figures/tikz/en_poisson_o1.txt};
        \addplot[color=red,dashed,mark=o,mark options={solid}] table [x index=3, y index=9] {figures/tikz/en_poisson_o1.txt};
        \addplot[color=blue,mark=square,solid,mark options={solid}] table [x index=3, y index=5] {figures/tikz/en_poisson_o2.txt};
        \addplot[color=blue,mark=square,dashed,mark options={solid}] table [x index=3, y index=9] {figures/tikz/en_poisson_o2.txt};
        \addplot[color=black!60!green,mark=square,solid,mark options={solid}] table [x index=3, y index=5] {figures/tikz/en_poisson_o3.txt};
        \addplot[color=black!60!green,mark=square,dashed,mark options={solid}] table [x index=3, y index=9] {figures/tikz/en_poisson_o3.txt};
        \addlegendentry{Initial ($p=1$)};
        \addlegendentry{Optimized ($p=1$)};
        \addlegendentry{Initial ($p=2$)};
        \addlegendentry{Optimized ($p=2$)};
        \addlegendentry{Initial ($p=3$)};
        \addlegendentry{Optimized ($p=3$)};
      \end{axis}
    \end{tikzpicture}
    }
    \caption{Load functional-based measure}
    \label{fig:energy_quad_l2}
  \end{subfigure}
  \newline
    \begin{subfigure}[b]{0.48\textwidth}
    \centering
    \resizebox{1.0\textwidth}{!}{%
    \begin{tikzpicture}
      \begin{axis}[
          legend columns=2,
          xlabel={$h$},
          ylabel={$||u_h-u^*||_{L_2}$},
          legend pos=north east,
          x dir=reverse,
          xmode=log,
          ymode=log,
          log basis x=10,
          log basis y=10,
          ymax=1e2,
          ymin=1e-8,
          legend cell align={left},
          every axis plot/.append style={line width=1pt}
        ]
        \addplot[color=red,mark=o,solid] table [x index=3, y index=5] {figures/tikz/zz_poisson_o1.txt};
        \addplot[color=red,dashed,mark=o,mark options={solid}] table [x index=3, y index=9] {figures/tikz/zz_poisson_o1.txt};
        \addplot[color=blue,mark=square,solid,mark options={solid}] table [x index=3, y index=5] {figures/tikz/zz_poisson_o2.txt};
        \addplot[color=blue,mark=square,dashed,mark options={solid}] table [x index=3, y index=9]  {figures/tikz/zz_poisson_o2.txt};
        \addplot[color=black!60!green,mark=square,solid,mark options={solid}] table [x index=3, y index=5] {figures/tikz/zz_poisson_o3.txt};
        \addplot[color=black!60!green,mark=square,dashed,mark options={solid}] table [x index=3, y index=9]  {figures/tikz/zz_poisson_o3.txt};
        \addlegendentry{Initial ($p=1$)};
        \addlegendentry{Optimized ($p=1$)};
        \addlegendentry{Initial ($p=2$)};
        \addlegendentry{Optimized ($p=2$)};
        \addlegendentry{Initial ($p=3$)};
        \addlegendentry{Optimized ($p=3$)};
      \end{axis}
    \end{tikzpicture}
    }
    \caption{Gradient-continuity measure}
    \label{fig:zz_quad_l2}
  \end{subfigure}
\caption{Convergence of the solution error for meshes adapted
         by different measures introduced in Section \ref{sec::PDEperformanceMeasure}.}
\label{fig:poisson_convergence}
\end{figure}
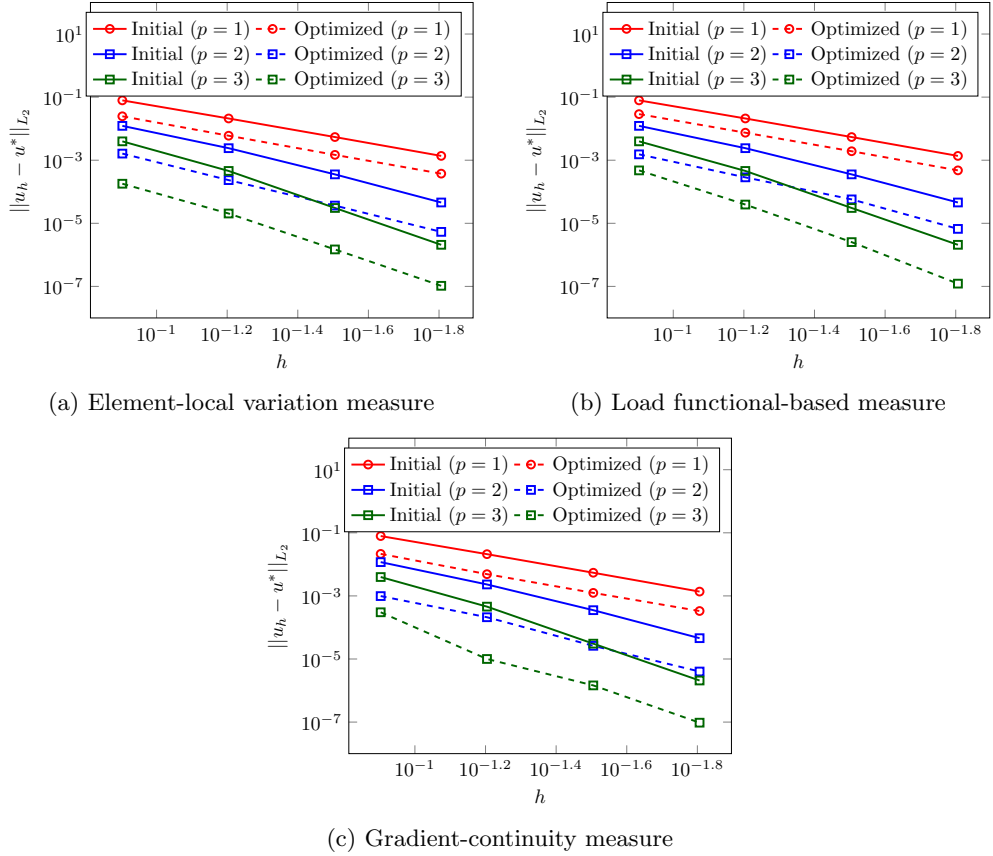

Figure \ref{fig_poisson_circular}(b) shows the $16\times16$ quadratic mesh optimized using the three error measures. We observe that in each case, the mesh size is reduced around the region of high-gradient as prescribed by the error measures, at the expense of resolution in less critical areas of the domain.
Figure \ref{fig:poisson_convergence} shows that the proposed methodology recovers the expected convergence rate of $\mathcal{O}(h^{p+1})$ by roughly reducing the error by a factor of 10 relative to the initial uniform mesh. Here, $h$ is the element size, and $p$ is the polynomial order used to represent both the mesh nodes and the finite element solution. Note that in these results, the penalization weight $\alpha$, filter radius $\delta$ and the number of solver iterations have been tuned to achieve the best possible performance. This is reasonable as the goal of this experiment is to demonstrate the effectiveness of the proposed framework. In future work, we plan to develop techniques to automate the selection of these parameters, for example, similar to penalization weight for mesh fitting in Figure \cite{barrera2023high}.

\begin{figure}[htbp]
\begin{center}
$\begin{array}{ccc}
\includegraphics[height=0.3\textwidth]{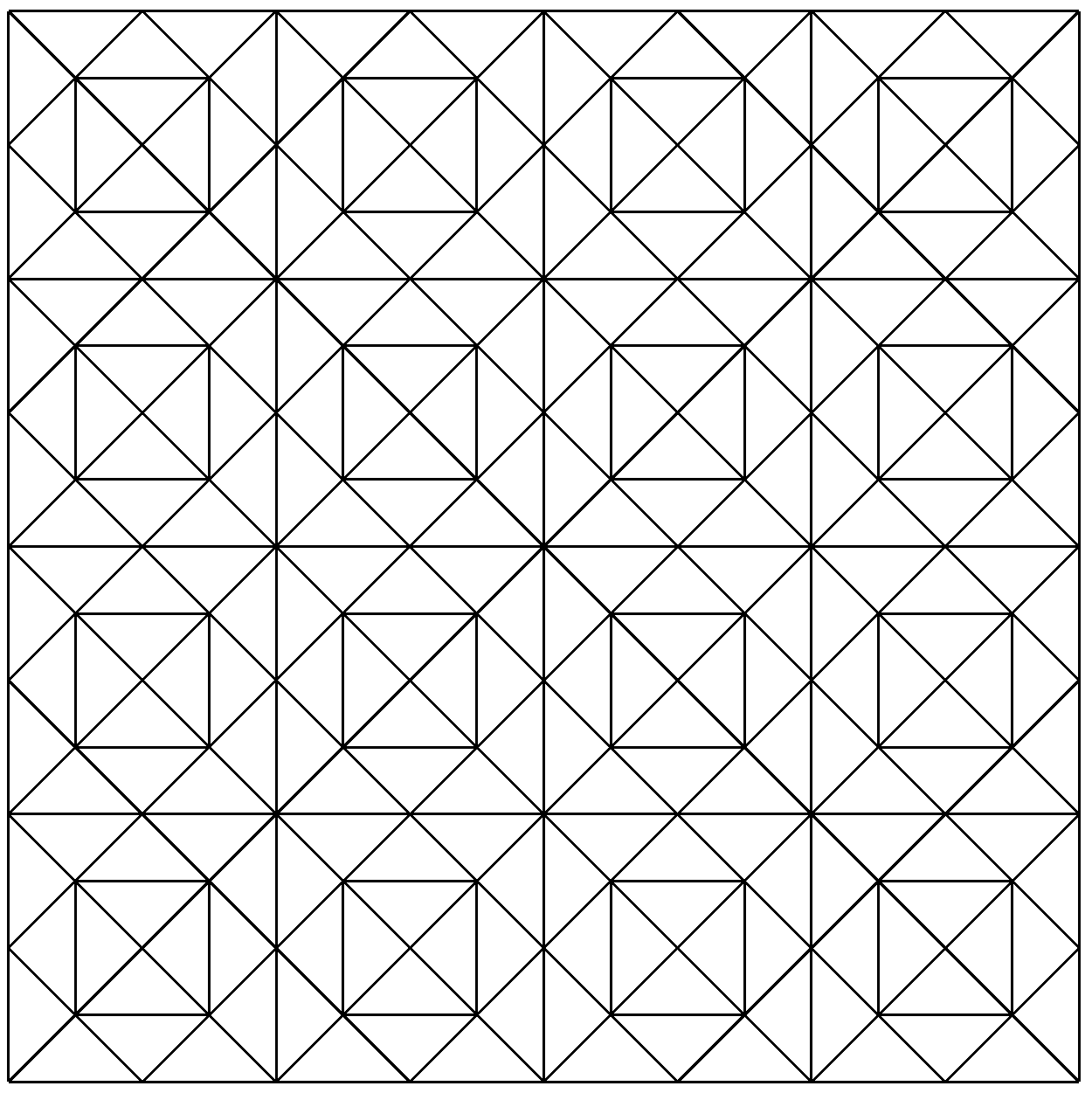} &
\includegraphics[height=0.3\textwidth]{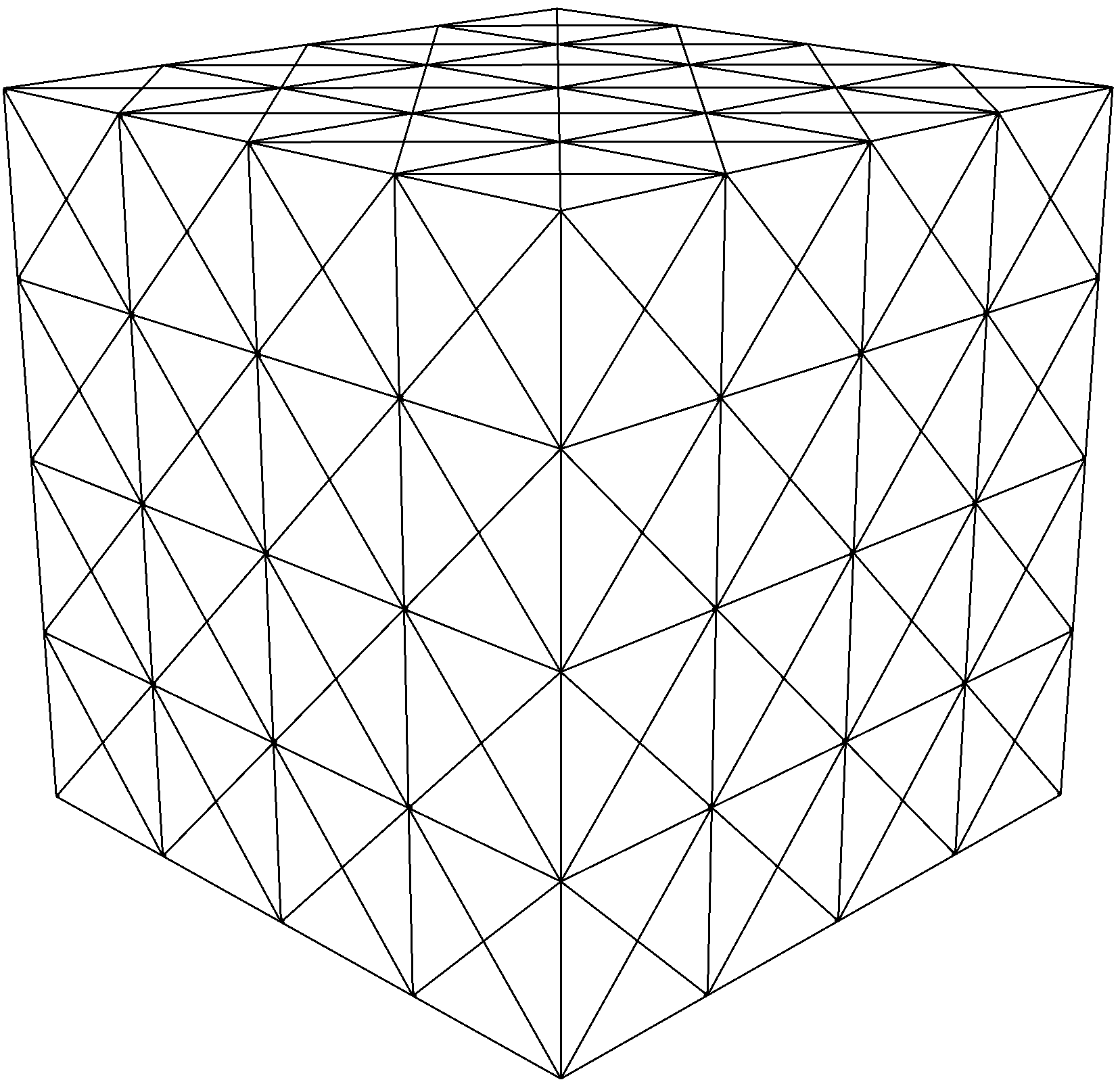} &
\includegraphics[height=0.3\textwidth]{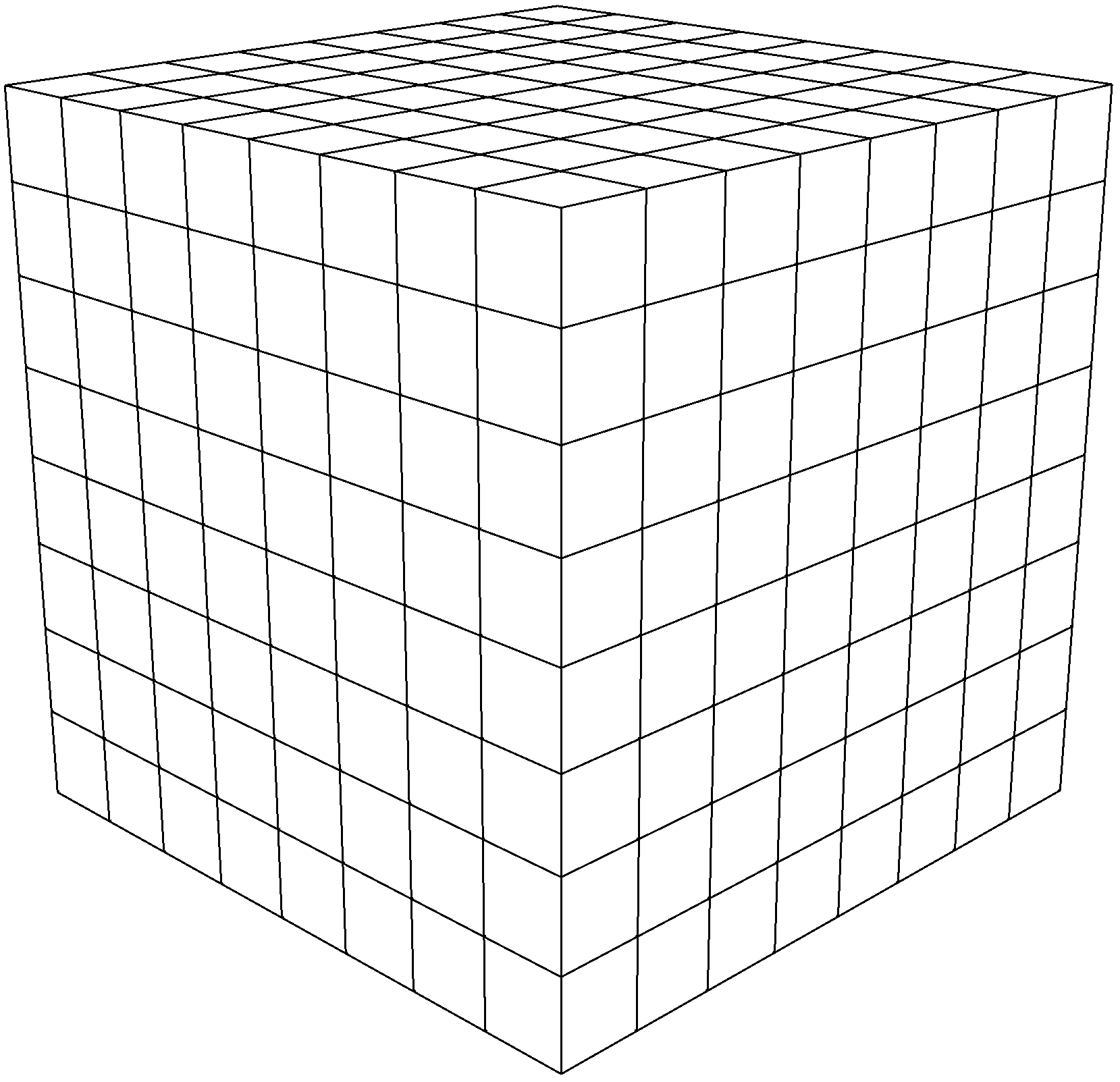} \\
\vspace{-2mm} \\
\textrm{(a)} & \textrm{(c)} & \textrm{(e)} \\
\includegraphics[height=0.3\textwidth]{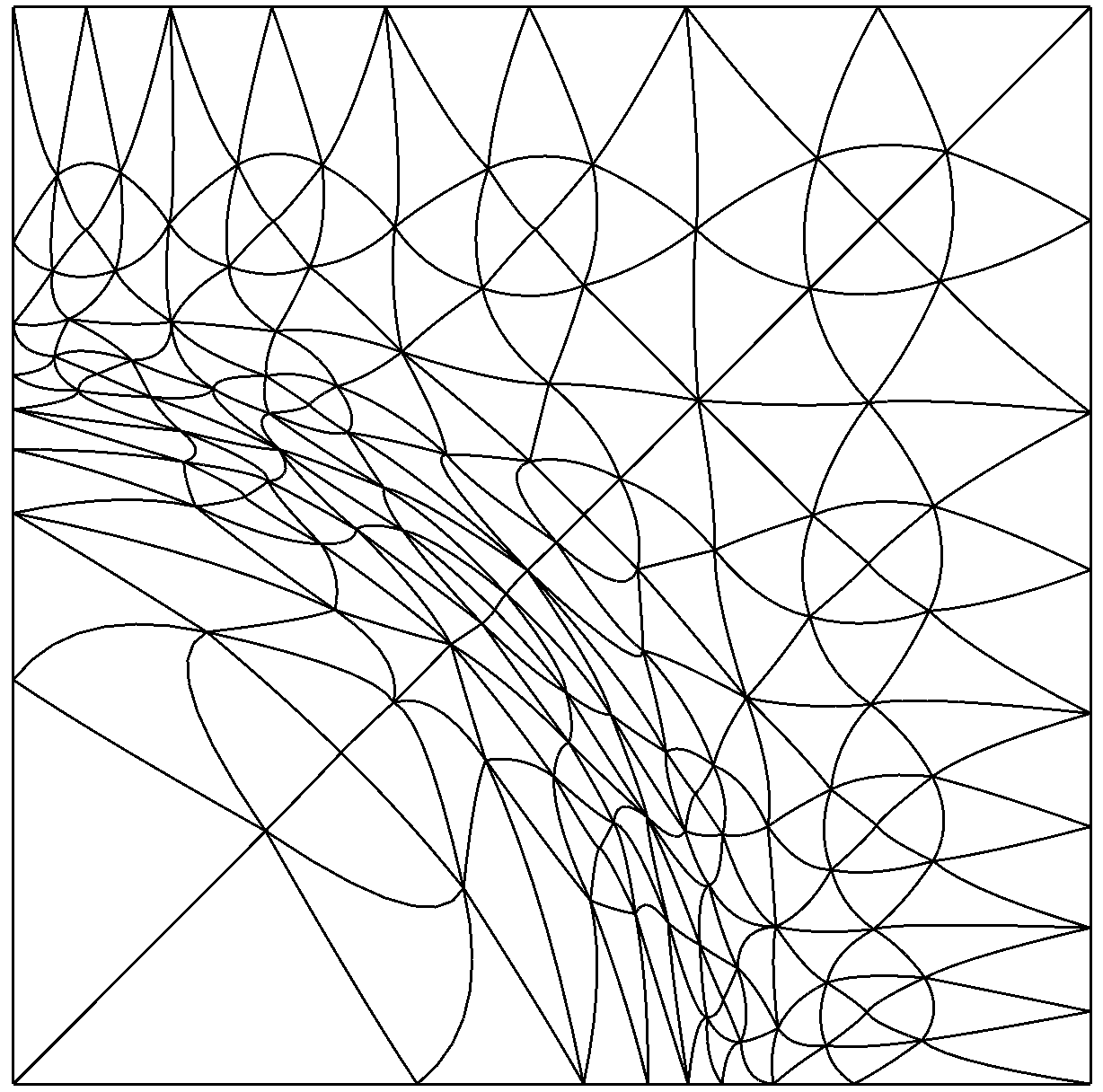} &
\includegraphics[height=0.3\textwidth]{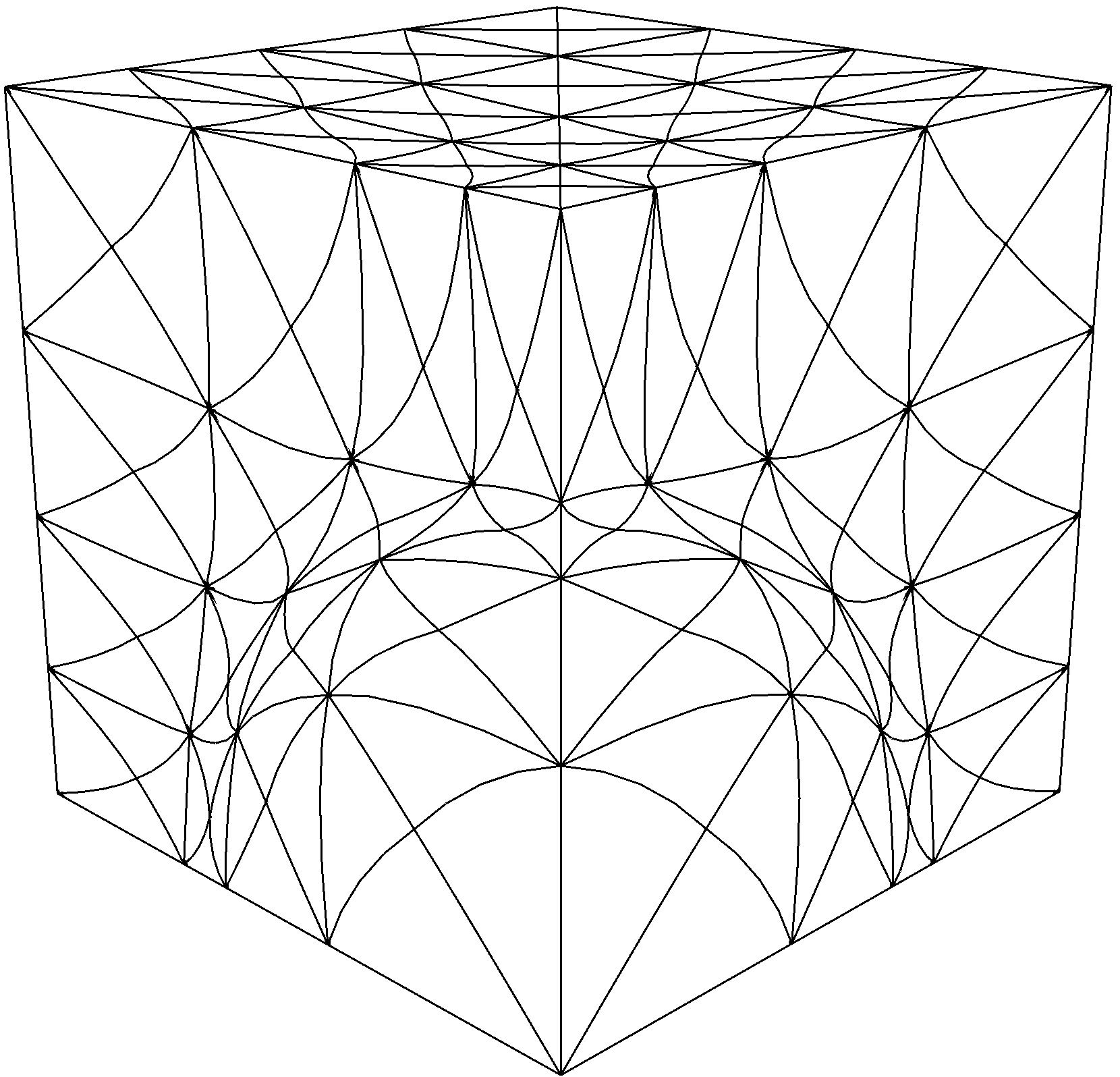} &
\includegraphics[height=0.3\textwidth]{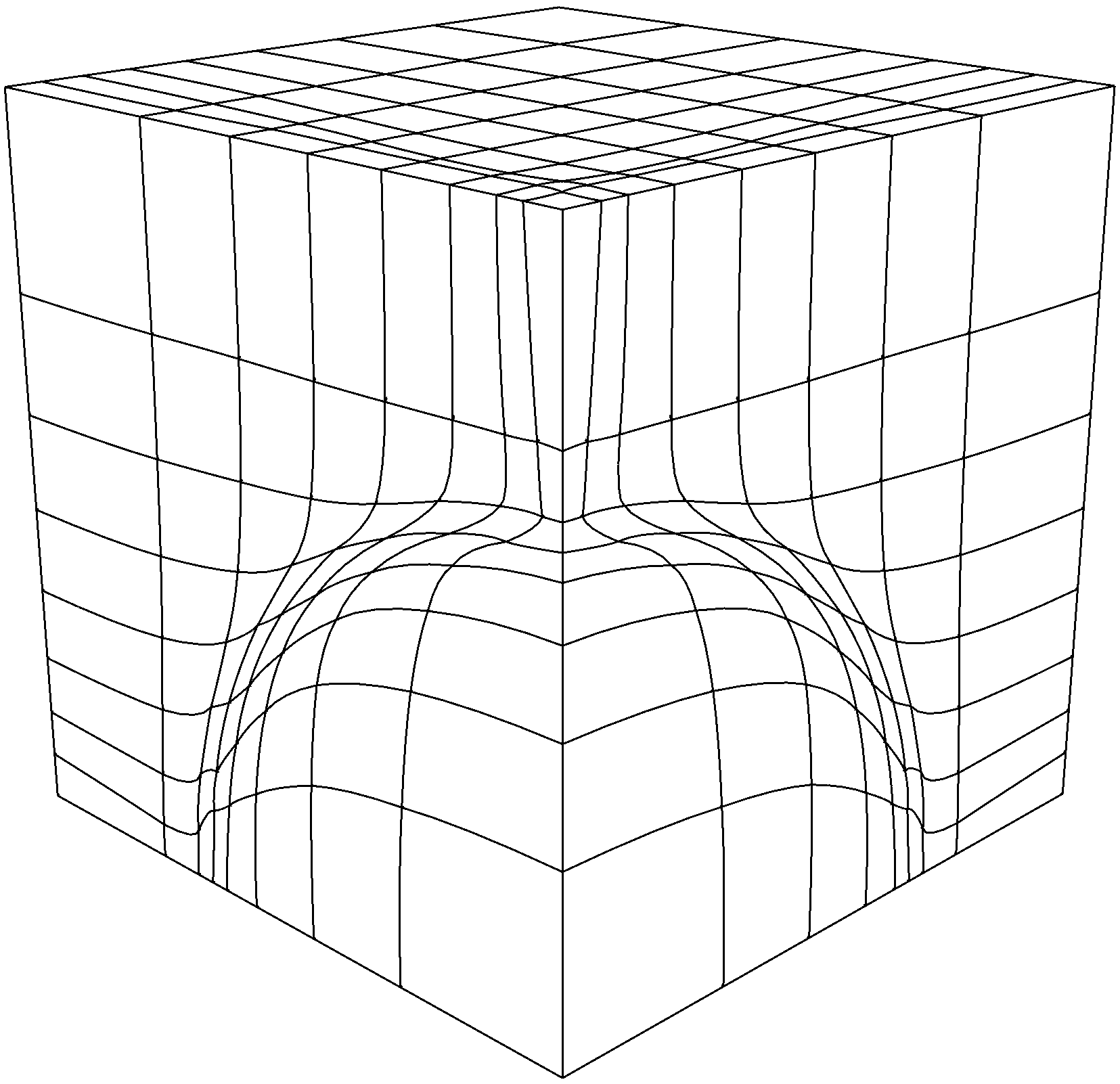} \\
\vspace{-2mm} \\
\textrm{(b)} & \textrm{(d)} & \textrm{(f)} \\
\end{array}$
\end{center}
% \vspace{-7mm}
\caption{Initial and optimized meshes with different element types (triangles, tets, hexes), adapted using the element-local variation measure.}
\label{fig_poisson_circular_element_types}
\end{figure}

Next, we apply the solution-average based error measure to the Poisson problem using quadratic triangular, hexahedral, and tetrahedral meshes, demonstrating the method’s generalization across different element types and spatial dimensions. The initial triangular mesh in Figure \ref{fig_poisson_circular_element_types}(a) is generated using a Cartesian aligned $4 \times 4$ quad mesh. Each quad element is split into four triangles by connecting the element center to the vertices. Each triangle is further refined into four triangles by connecting the midpoint of its three edges. Similarly, the tetrahedral mesh is generated by taking a $4 \times 4 \times 4$ hexahedral mesh and splitting each hex into 24 tetrahedra. This decomposition is achieved by dividing each of the six hexahedral faces into four triangles and connecting the resulting triangles to the centroid of the hex. The initial tetrahedral mesh is shown in Figure \ref{fig_poisson_circular_element_types}(c). The initial hexahedral mesh is an $8 \times 8 \times 8$ Cartesian-aligned mesh, and it is shown in Figure \ref{fig_poisson_circular_element_types}(e). In each case, the mesh quality measure is based on a shape metric with a target matrix set such that the ideal element is a cube for hex mesh and an equilateral simplex for the simplicial meshes. For all cases, we observe that the proposed PDE-constrained approach results in a mesh with element size reduced along the curve featuring a high gradient in the solution.

\begin{figure}[htbp]
\begin{center}
$\begin{array}{ccc}
\includegraphics[height=0.3\textwidth]{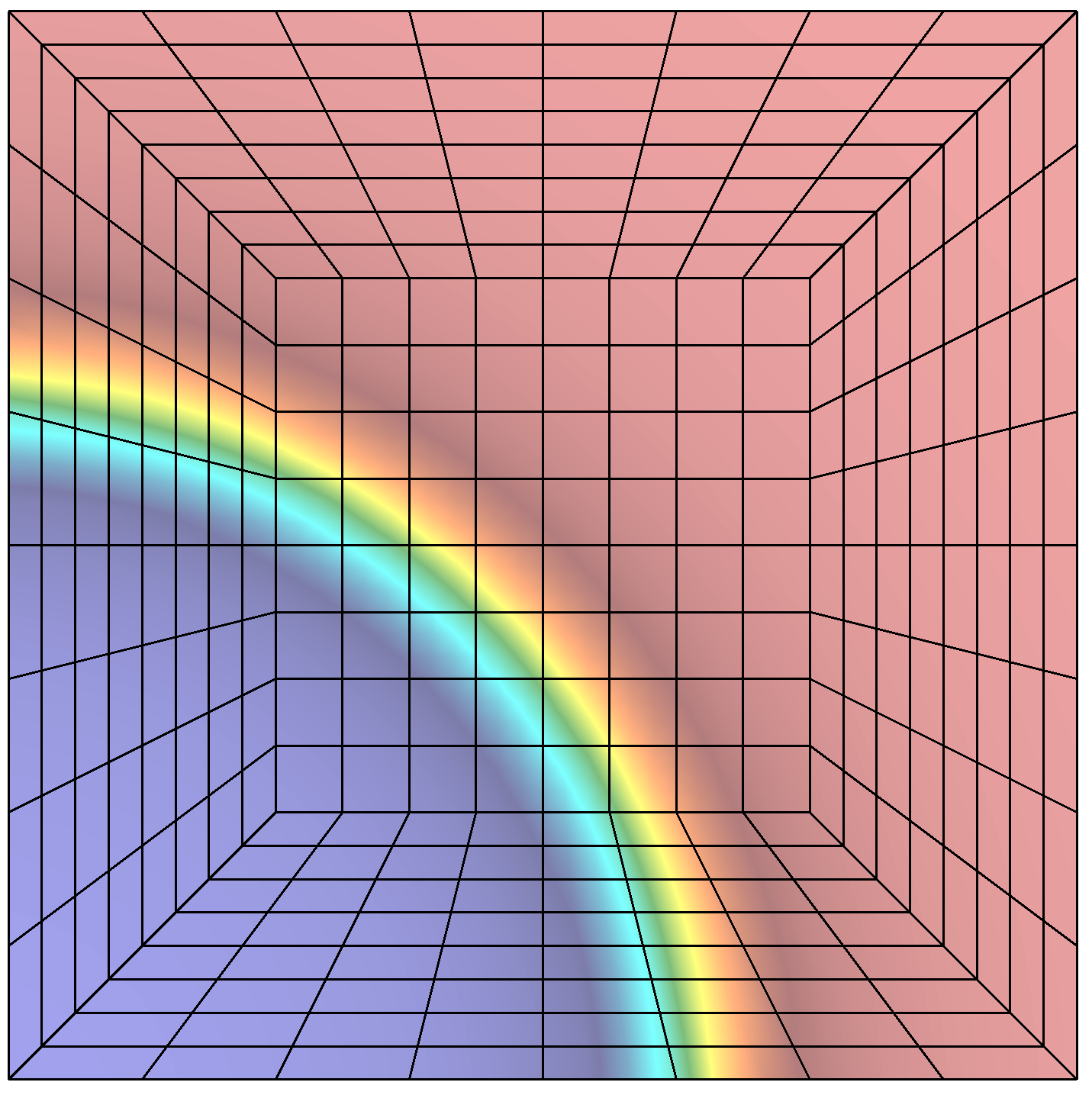} &
\includegraphics[height=0.3\textwidth]{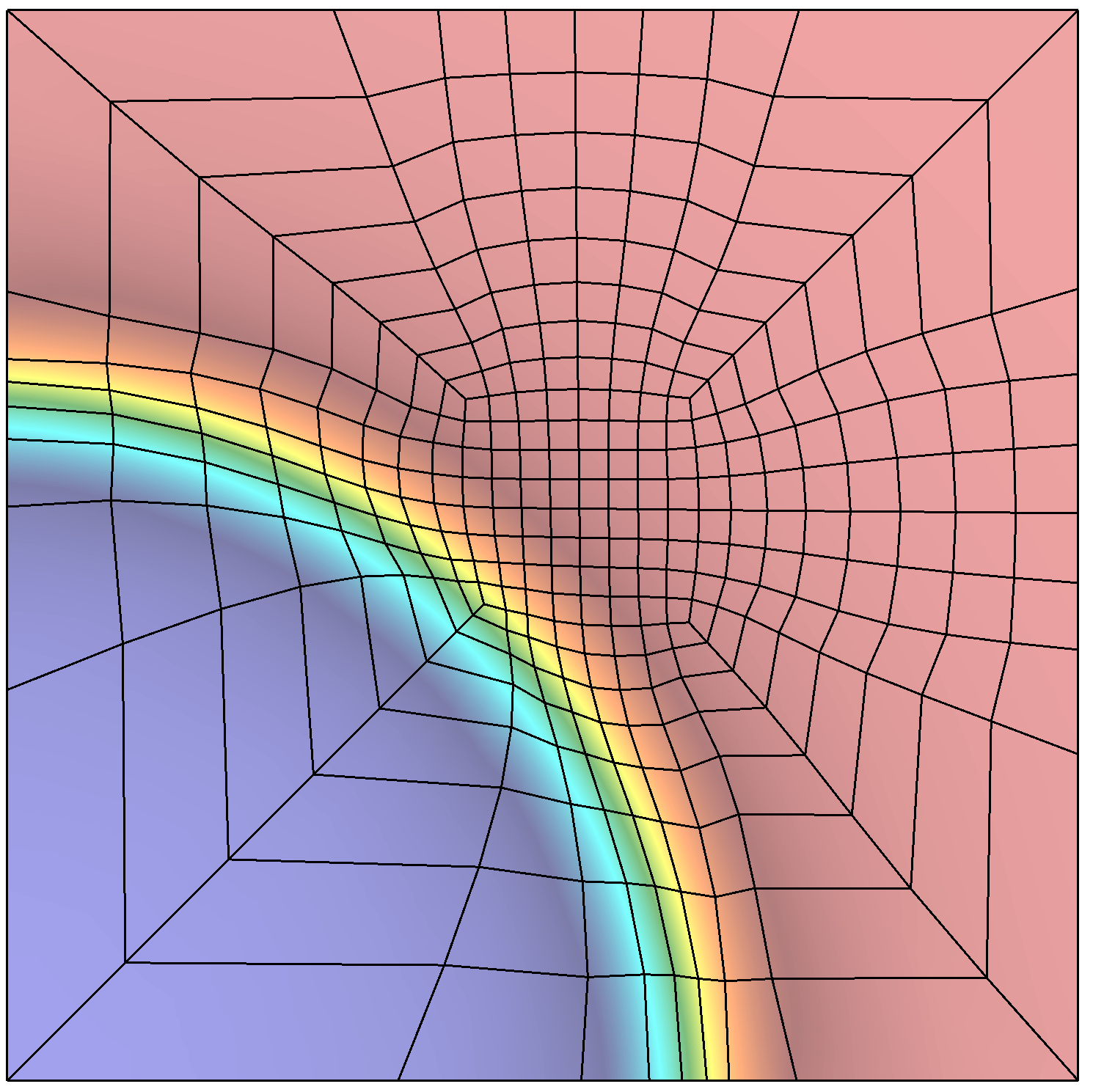} &
\includegraphics[height=0.3\textwidth]{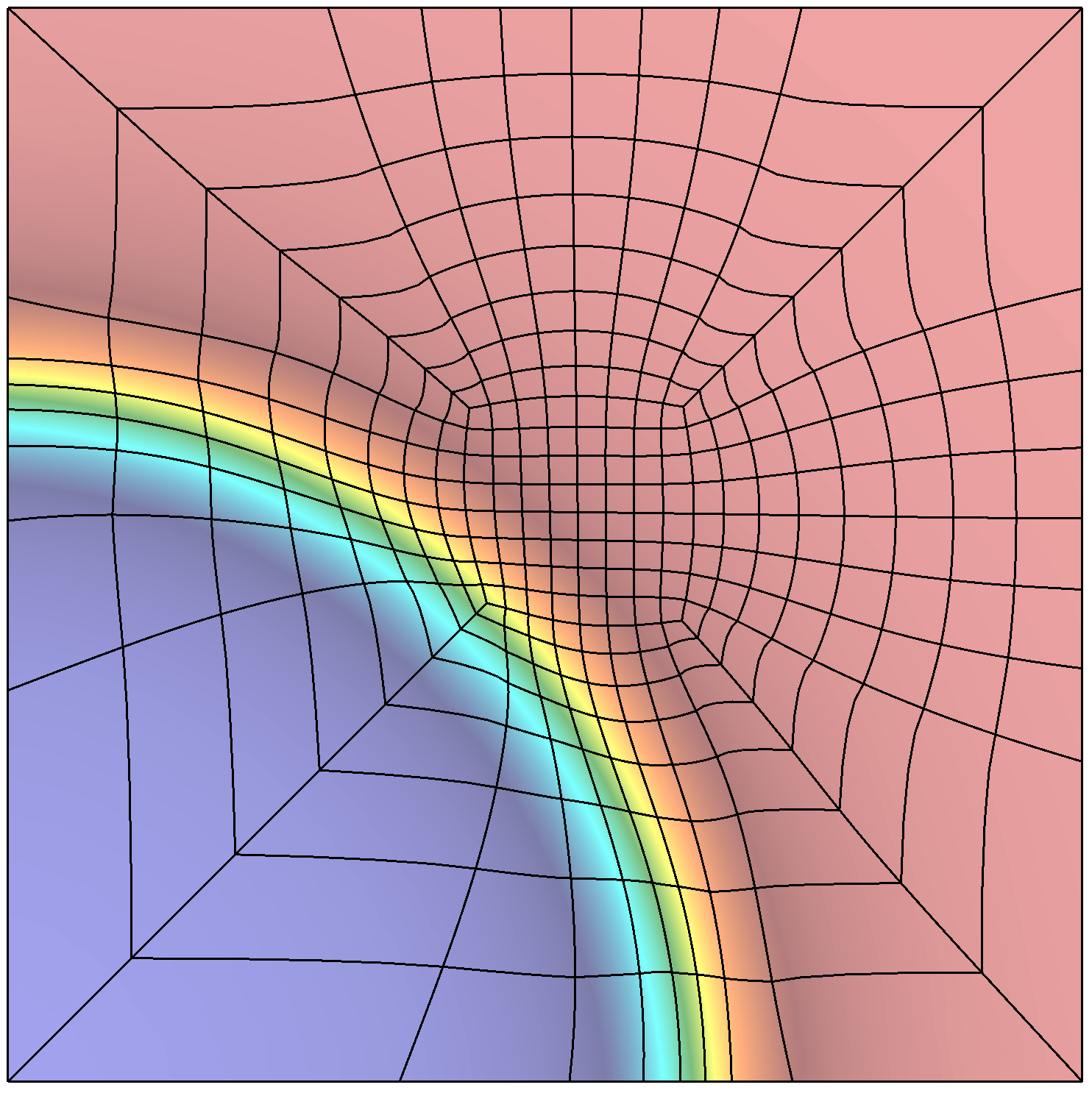} \\
\vspace{-2mm} \\
\textrm{(a) Initial mesh} & \textrm{(b) Optimized linear mesh} & \textrm{(c) Optimized quadratic mesh}
\end{array}$
\end{center}
% \vspace{-7mm}
\caption{Adapting an unstructured mesh using the element-local variation measure with the proposed framework. PDE solution show on the (a) Initial mesh, (b) optimized linear mesh ($p=1$), and (c) optimized quadratic mesh ($p=2$). In each case, the PDE solution is based on a second-order finite element space ($p_u=2$).}
\label{fig_poisson_unstructured}
\end{figure}

Finally, Figure \ref{fig_poisson_unstructured} shows the effectiveness of the proposed framework for an unstructured mesh where the polynomial order used to represent the mesh is different from that of the solution. This example utilizes the load functional-based measure from Section \ref{sec::compliance} in $\mcf_P$. We observe that the optimized linear and quadratic meshes have similar element sizes throughout the domain, but the quadratic mesh is able to better align with the region of high gradients due to increased flexibility from additional degrees of freedom per element. As a result, the quadratic optimized mesh has lower error in the solution, measured in terms of both the $L^2$ norm and the $H^1$ semi-norm.

\subsection{Mesh alignment with a vector field}
\begin{figure}[htbp]
\centering
\begin{minipage}{0.3\textwidth}
  \centering
  \includegraphics[width=\textwidth]{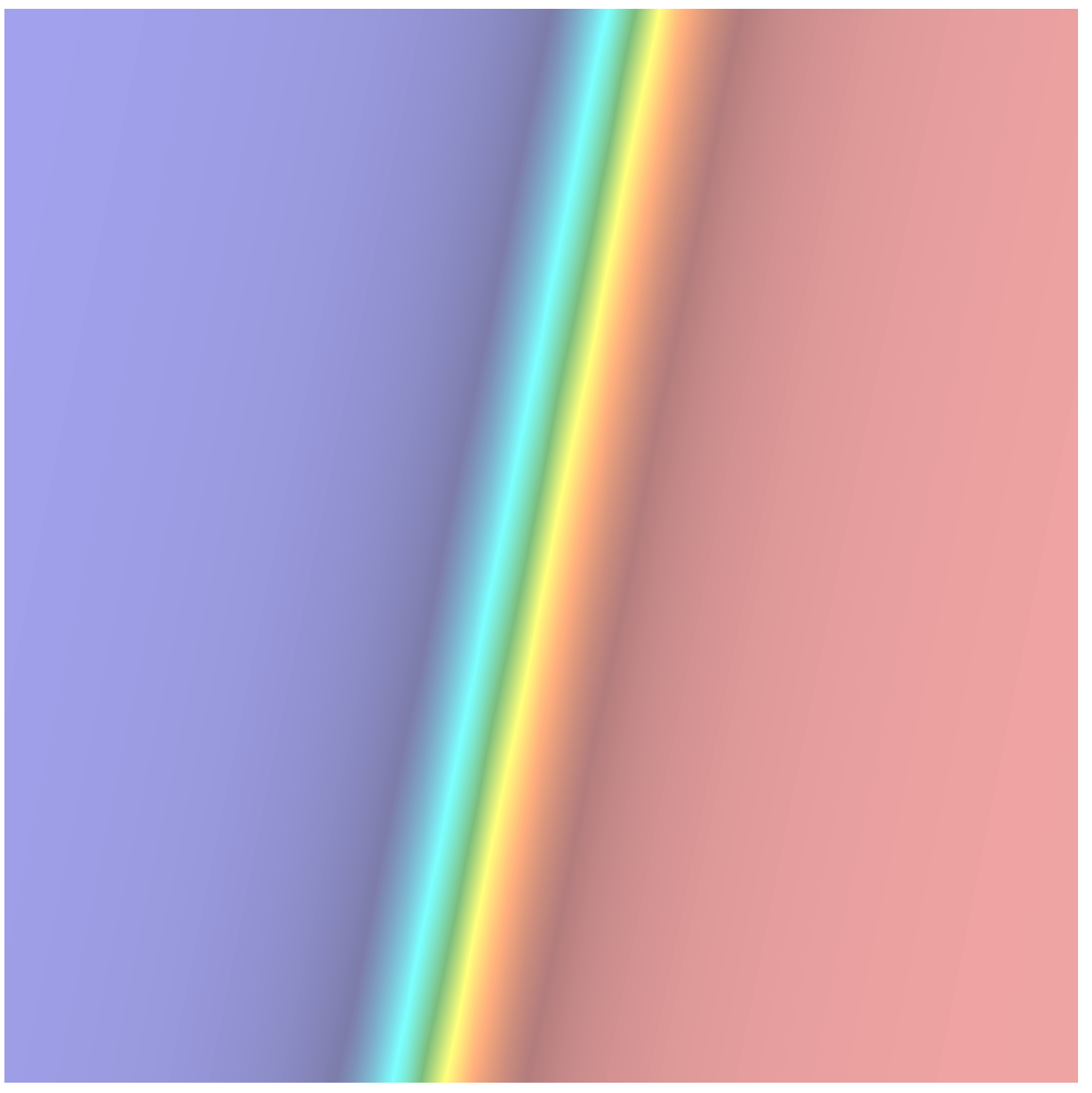}
(a)
\end{minipage}
\hspace{2em}
\begin{minipage}{0.3\textwidth}
  \centering
  \includegraphics[width=\textwidth]{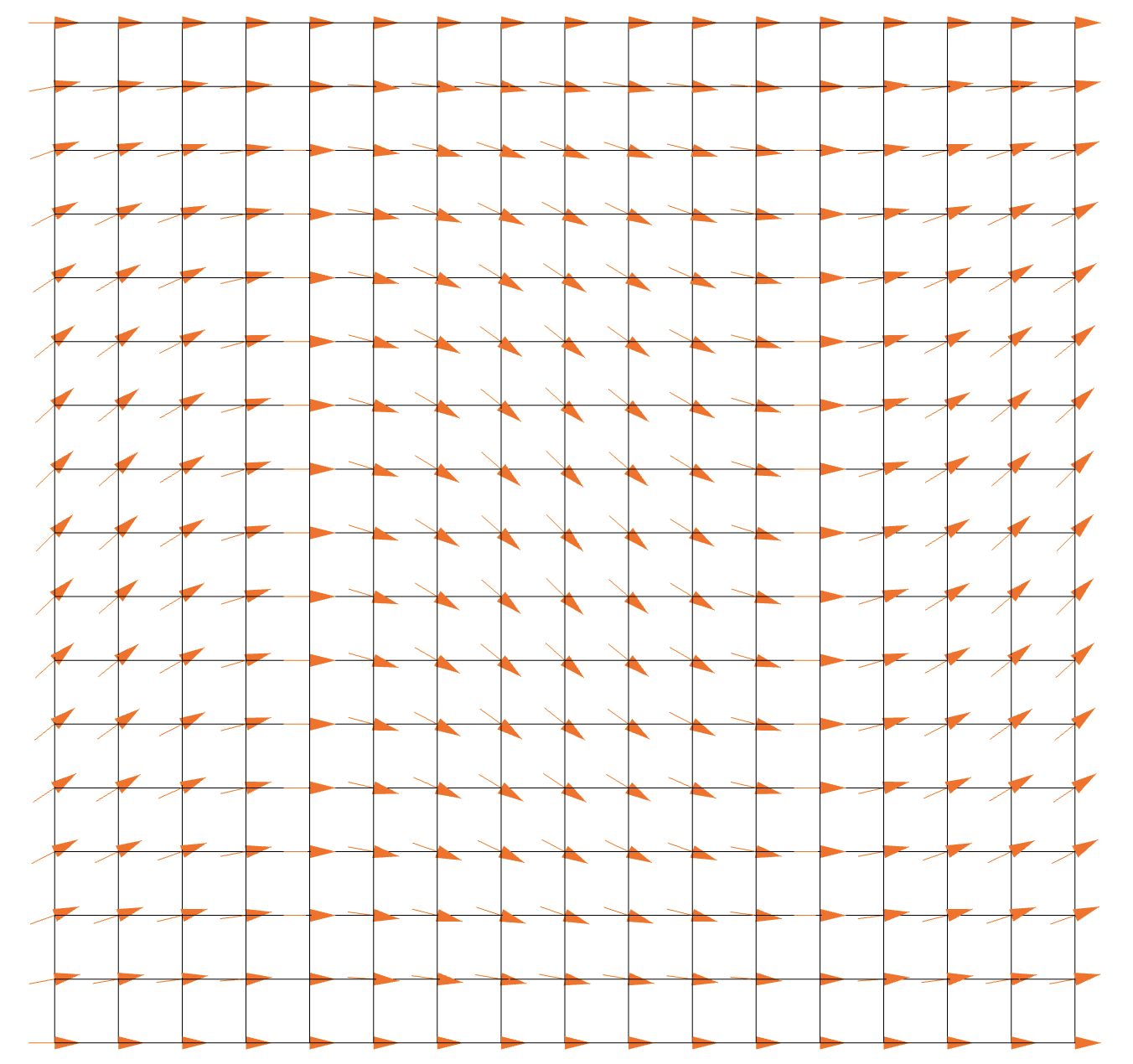}
  (b)
\end{minipage}
\begin{minipage}{0.3\textwidth}
  \centering
  \includegraphics[width=\textwidth]{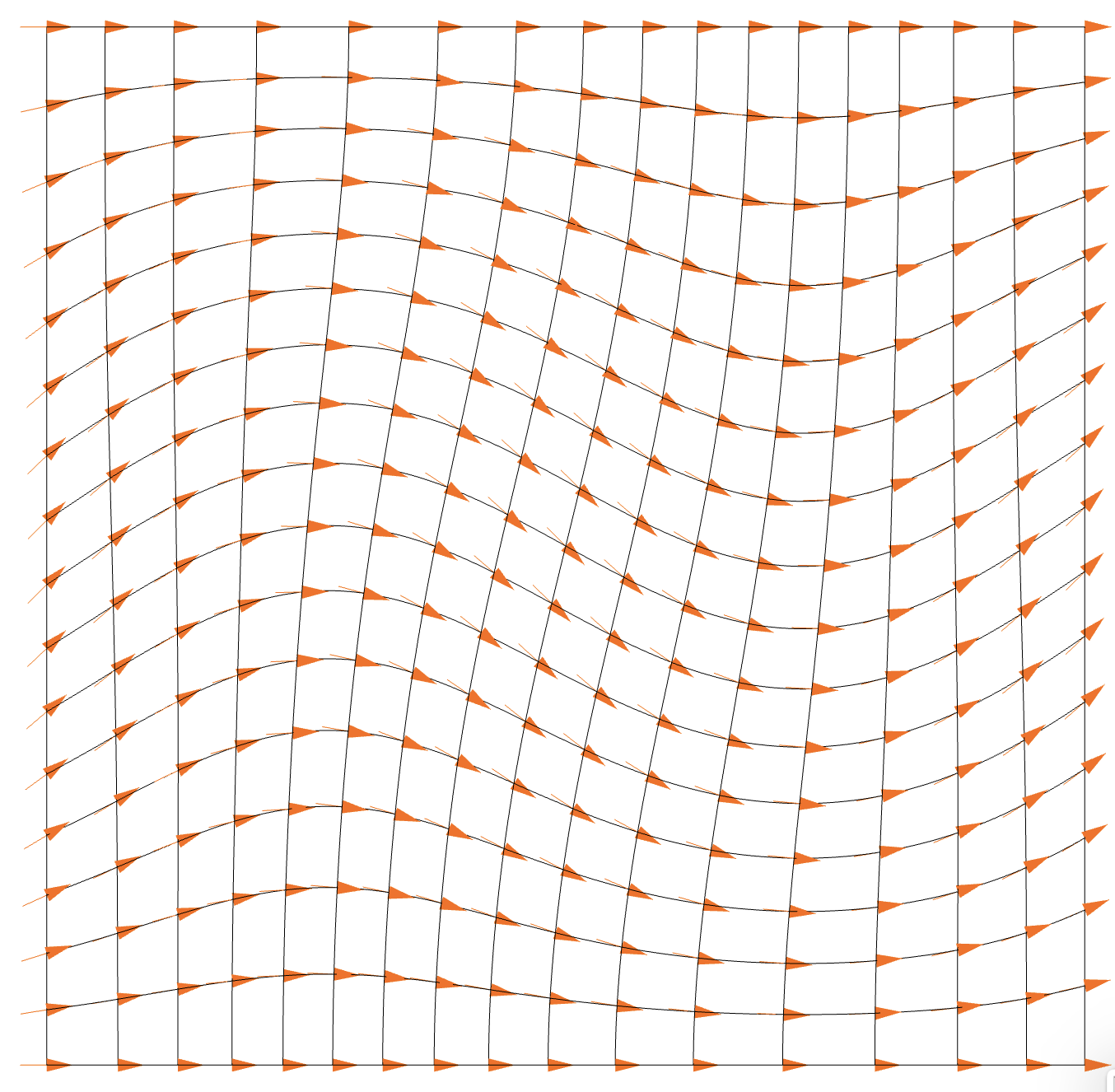}
  (c) \(\alpha=0\)
\end{minipage}
\newline
\hspace{1em}
\begin{minipage}{0.3\textwidth}
  \centering
  \includegraphics[width=\textwidth]{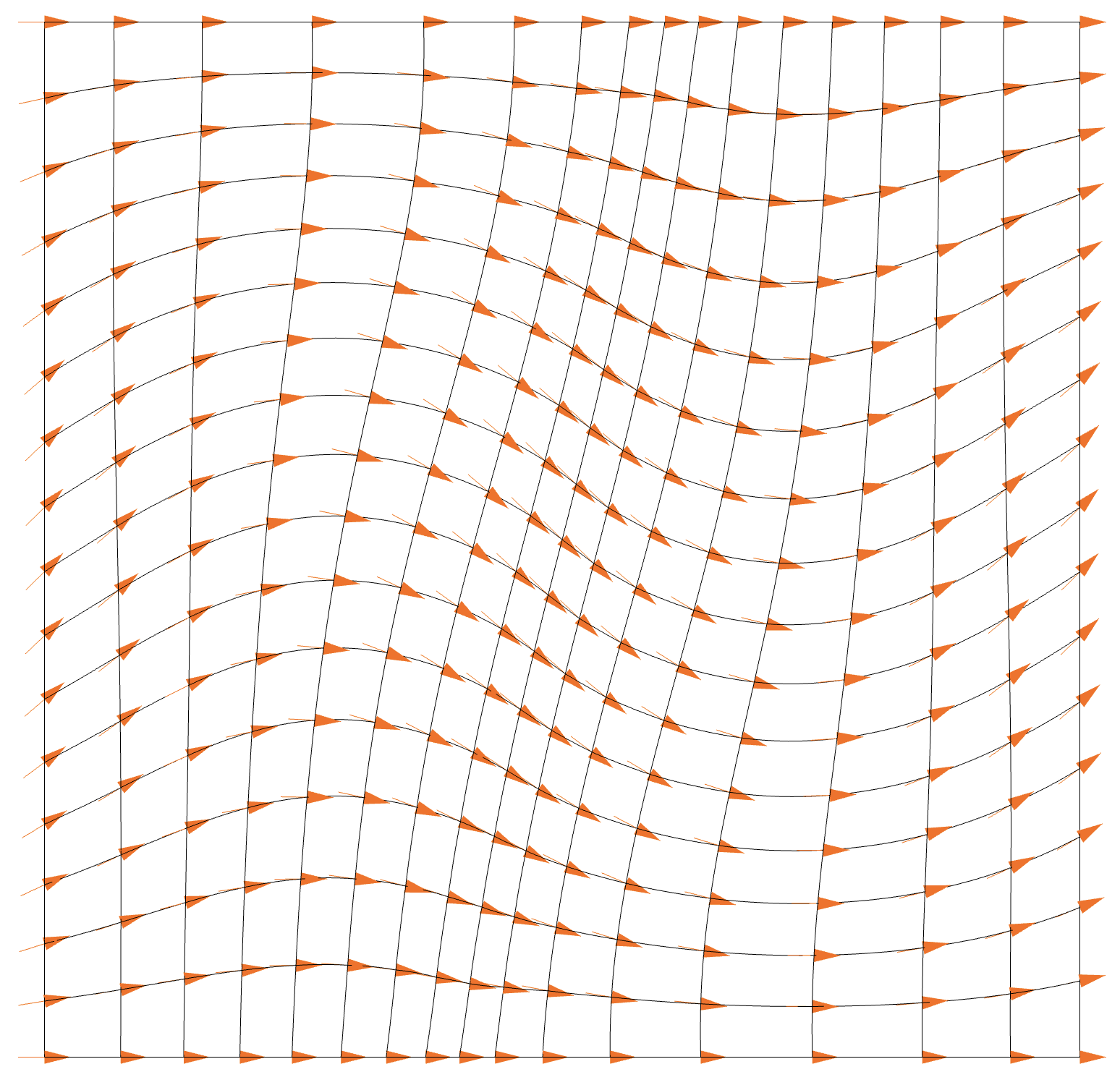}
  (d) \(\alpha=4\cdot 10^3\)
\end{minipage}
\hspace{1em}
\begin{minipage}{0.3\textwidth}
  \centering
  \includegraphics[width=\textwidth]{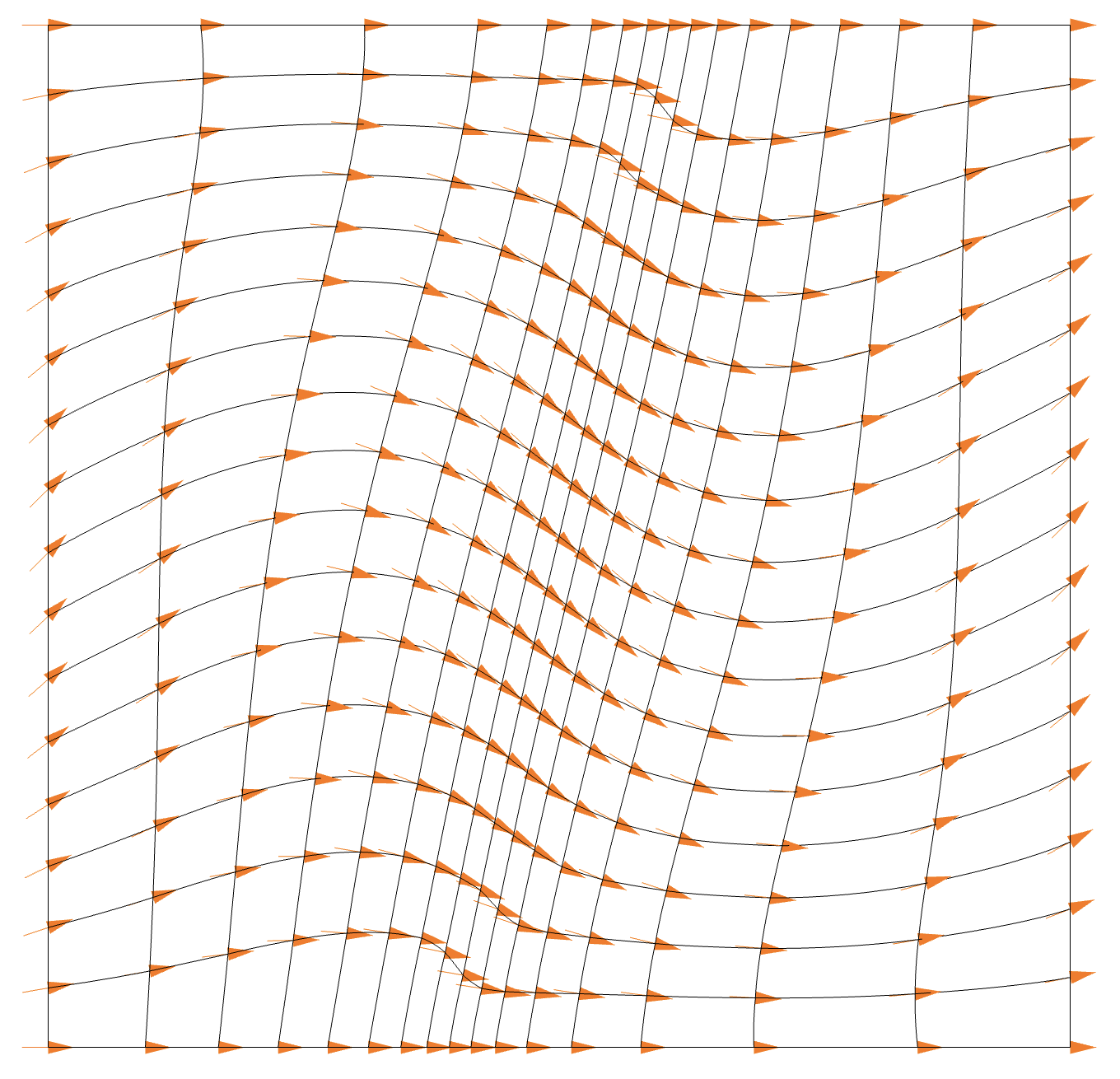}
  (e) \(\alpha=4\cdot 10^4\)
\end{minipage}
\hspace{1em}
\begin{minipage}{0.3\textwidth}
  \centering
  \includegraphics[width=\textwidth]{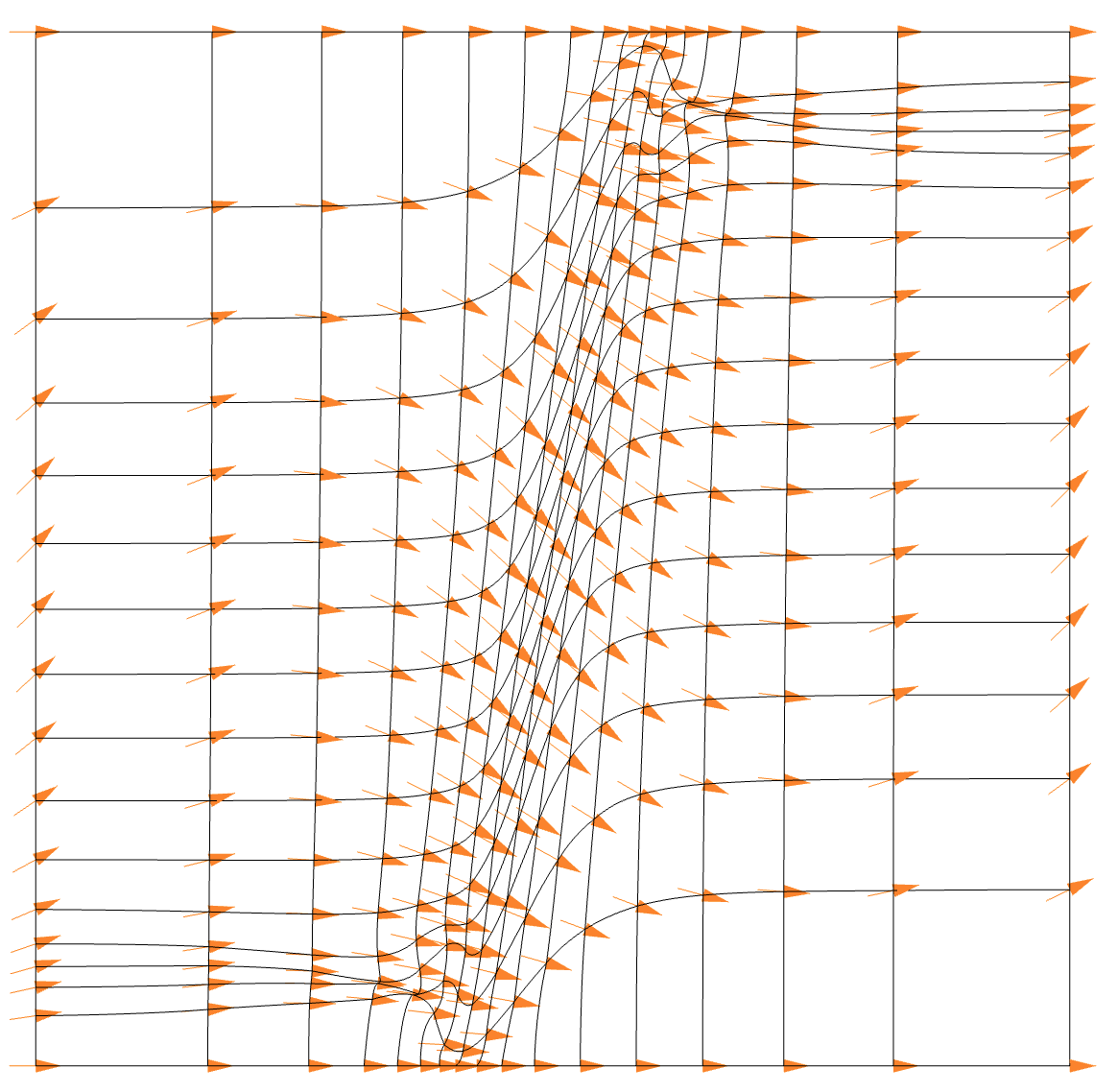}
  (f) \(\alpha=4\cdot 10^6\)
\end{minipage}
\caption{Comparison of meshes obtained with different penalization weight for the element-local solution variation measure. (a) Exact solution of the Poisson problem \eqref{eq:poisson_inclined_solution}. (b) Initial mesh with the target element orientation encoded by the target matrix $W$ in \eqref{eq:inclined_target}. Optimized meshes obtained with (c) $\alpha=0$, (d) $\alpha=4\cdot 10^3$, (e) $\alpha=4\cdot 10^4$, and (f) $\alpha=4\cdot 10^6$.}
\label{fig_poisson_inclined}
\end{figure}

The following example demonstrates the effectiveness of the proposed method on the Poisson equation, using a mesh quality measure that aligns the mesh with a prescribed vector field. The exact solution of the Poisson problem is chosen to mimic a shock wave in a unit-sized domain, $\Omega \in [0,1]^2$:
\begin{subequations}
    \begin{align}
        \label{eq:poisson_inclined_solution}
        -\nabla^2 u &= f, \\
u(\bx) &= \tan^{-1}\bigg(\alpha\big(x-0.5-0.2(y-0.5)\big)\bigg).
    \end{align}
\end{subequations}
Inhomogeneous Dirichlet conditions are prescribed on all boundaries. The exact solution of the Poisson problem is shown in Figure \ref{fig_poisson_inclined}(a). For the mesh quality measure, the orientation component of the target matrix is based on the vector field shown in Figure \ref{fig_poisson_inclined}(b), while the shape components encode unity aspect-ratio and $\pi/2$ skewness:
\begin{subequations}
    \label{eq:inclined_target}
    \begin{align}
        W &= \begin{bmatrix}
    \cos\theta & \sin\theta \\
    -\sin\theta & \cos\theta
    \end{bmatrix}, \\
    \theta(\bx) &= \pi y (1 - y) \cos(2\pi x).
    \end{align}
\end{subequations}
The mesh quality measure uses the \emph{shape+orientation} metric $\nu_{107,OS}$, and the penalization weight for the solution-average error measure (Section \ref{sec::perf_sol_local}) is set to $\alpha=0$, $4\cdot 10^3$, $4\cdot 10^4$, or $4\cdot 10^6$ in \eqref{eq::linelast}. The filter radius is held fixed at $\delta^2=0.005$.

The initial mesh is a $16 \times 16$ uniform-sized Cartesian-aligned quad mesh. Using $\alpha=0$ results in a mesh that aligns with the vector field while maintaining good element shape. When $\alpha=4\cdot 10^3$, there is a 4 times reduction in the error of the PDE solution as compared to the original mesh, and the mesh is still well aligned with the target vector field. Increasing the penalization weight to $\alpha=4\cdot 10^4$ results in a 10-fold reduction in the error at the expense of the mesh alignment with the vector field in the region of high gradients in the solution. Increasing the penalization weight further to $\alpha=4\cdot 10^6$ results in a severely degrading mesh quality, which results in a higher error compared to the mesh optimized with $\alpha=4\cdot 10^4$. This behavior is evident through a visual inspection of the optimized meshes in Figure \ref{fig_poisson_inclined}(c)-(e), and a quantitative comparison of $\mcf_P$, $\mcf_{\mu}$, and the $L_2$ error in the solution in Table \ref{tab:mesh_comparison_inclined}. This example demonstrates the trade-off nature of the proposed objective, which gives the user the flexibility to set the penalization weight based on the desired goal. Extension of the proposed method to support PDEs such that the user can adapt the mesh based on multiple discrete solution fields (e.g., velocity and pressure for the Navier-Stokes equations) will be presented in future works.

\begin{table}[htbp]
    \centering
    \begin{tabular}{|l|c|c|c|}
        \hline
        \textbf{Mesh} & $\mcf_p$ & $F_{\mu}$ & $||u_{h}-u^*||_{L_2}$ \\
        \hline
        Initial & 0.01 & 40.9 & 0.0034 \\
        \hline
        Optimized $(\alpha = 0)$ & 0.008 & 23 & 0.0015 \\
        \hline
        Optimized $(\alpha = 4 \cdot 10^3)$ & 0.004 & 29.9 & 0.0007 \\
        \hline
        Optimized $(\alpha = 4 \cdot 10^4)$ & 0.002 & 60.1 & 0.00025 \\
        \hline
        Optimized $(\alpha = 4 \cdot 10^6)$ & 0.0014 & 393.7 & 0.00032 \\
        \hline
    \end{tabular}
    \caption{Comparison of $\mcf_p$, $F_{\mu}$, and $L_2$ errors for initial and optimized meshes with different penalization weight for the element-local solution variation measure.}
    \label{tab:mesh_comparison_inclined}
\end{table}

\subsection{Linear elasticity examples}

\begin{figure}[htbp]
\begin{center}
$\begin{array}{c}
\includegraphics[width=0.8\textwidth]{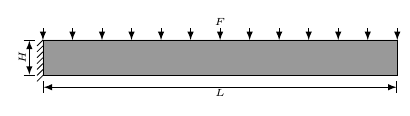}
\end{array}$
\end{center}
% \vspace{-7mm}
\caption{Two dimensional cantilever with boundary conditions.}
\label{fig_beam_sketch}
\end{figure}

\begin{figure}[htbp]
\begin{center}
$\begin{array}{c}
\includegraphics[width=0.4\textwidth]{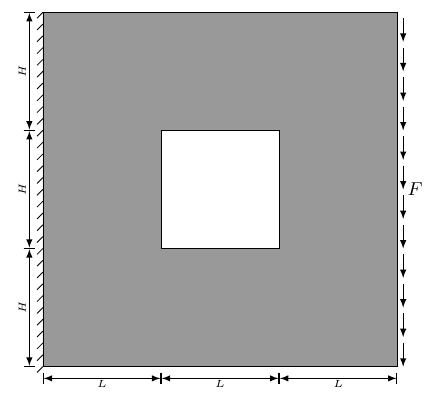}
\end{array}$
\end{center}
% \vspace{-7mm}
\caption{Sketch of two dimensional shear wall with a square hole and boundary conditions.}
\label{fig_platehole_sketch}
\end{figure}

The final two examples include a linear elastic beam subject to a uniform distributed load as illustrated in Figure~\ref{fig_beam_sketch}, and a square shear wall with a hole shown in Figure~\ref{fig_platehole_sketch}. These two examples demonstrate the applicability of the presented methodology to PDEs with vector-valued solutions. In both examples, we minimize the load functional defined in Section~\ref{sec::compliance} to obtain the best approximation among all discrete solutions over meshes with fixed topology.

\begin{figure}[htbp]
\begin{center}
$\begin{array}{c}
\includegraphics[width=0.8\textwidth]{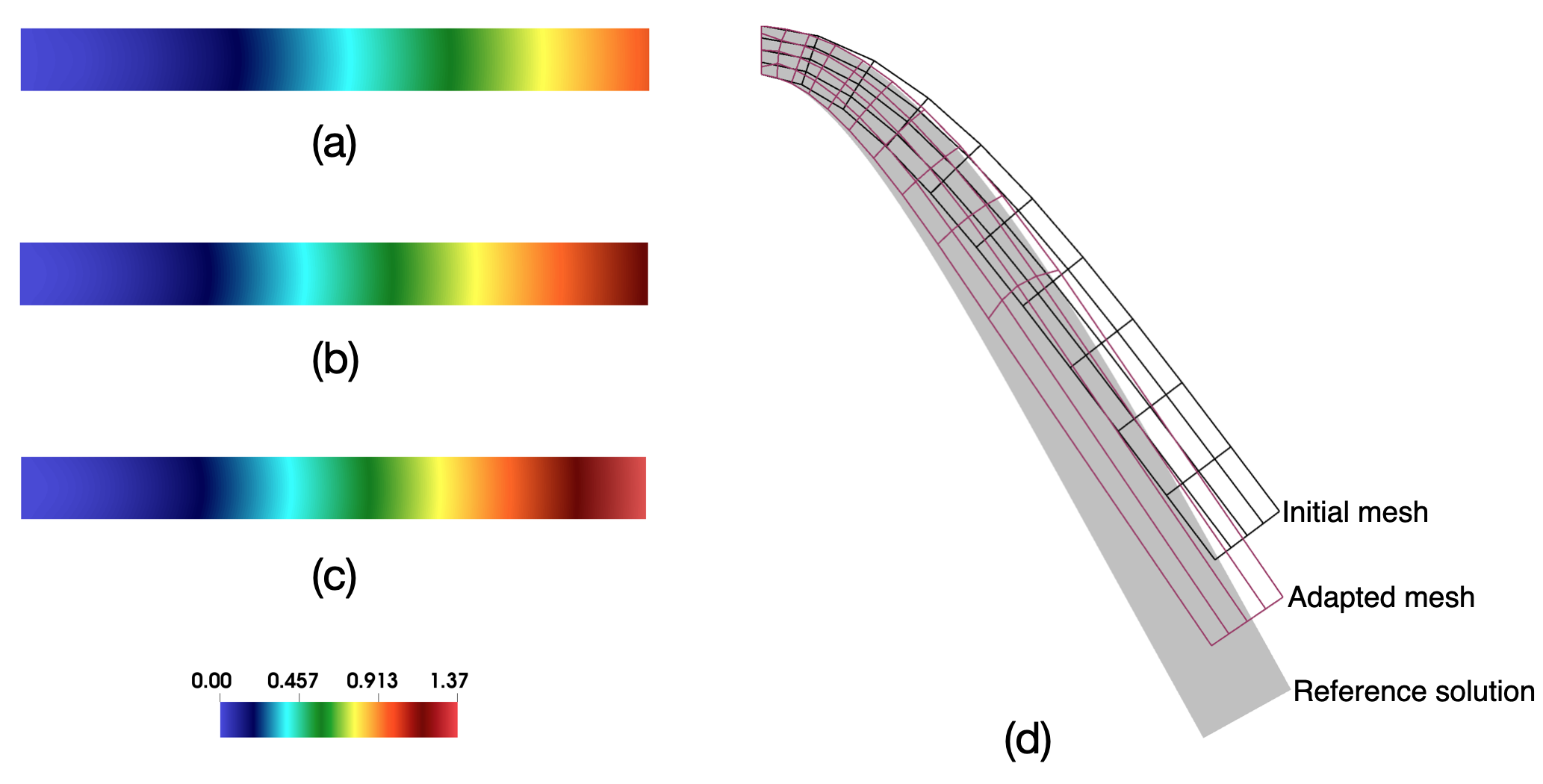}
\end{array}$
\end{center}
% \vspace{-7mm}
\caption{Displacement comparison of the initial and adapted linear meshes with the reference solution for the two dimensional cantilever beam. Displacement magnitude is shown for the (a) initial mesh, (b) adapted mesh, and (c) the reference solution. (d) The optimized mesh captures the physical response with higher fidelity compared to the initial mesh.}
\label{fig_beam_disp_mag_p1}
\end{figure}

\begin{figure}[htbp]
\begin{center}
$\begin{array}{cc}
\includegraphics[width=0.3\textwidth]
{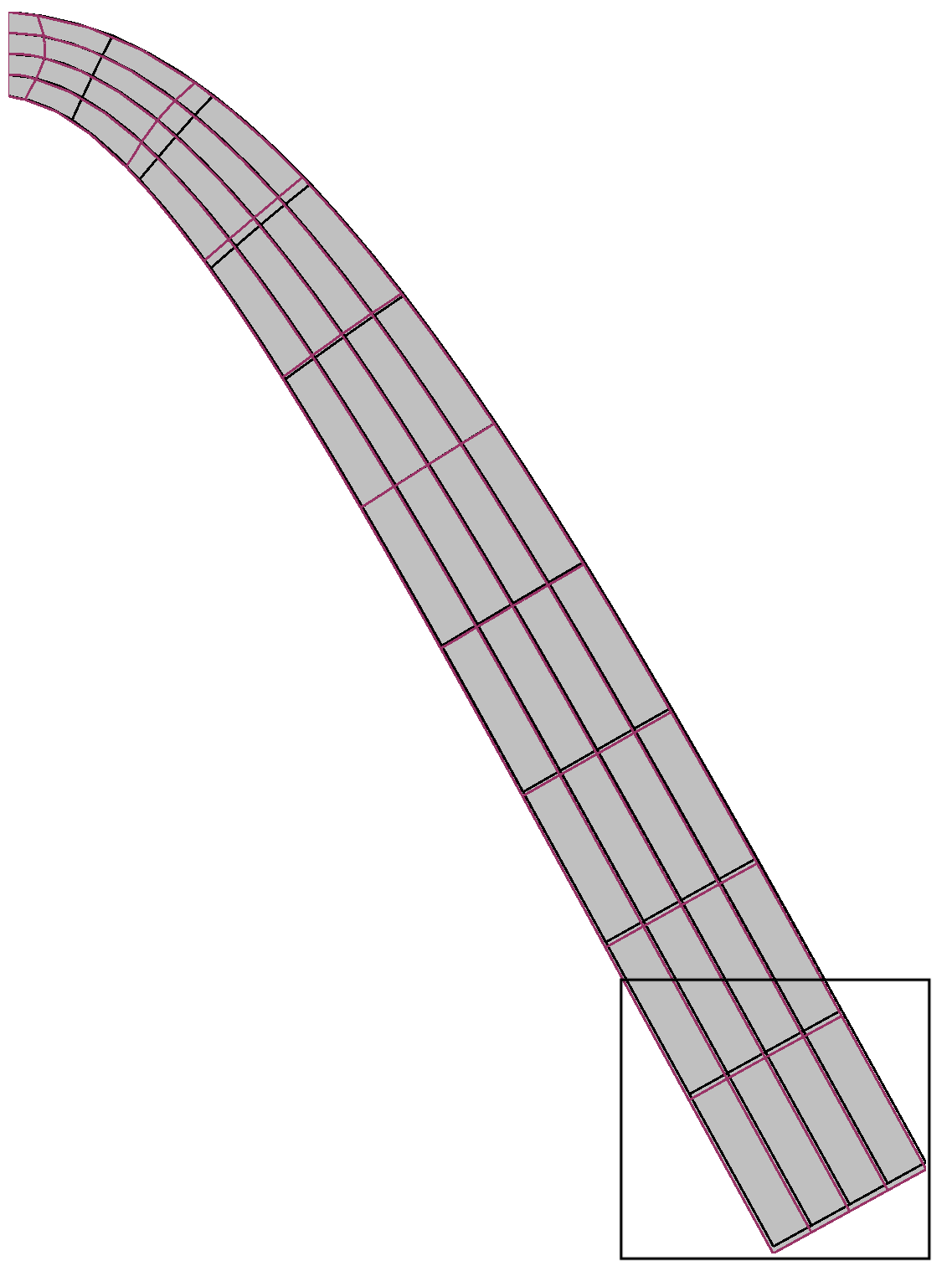} &
\includegraphics[width=0.4\textwidth]{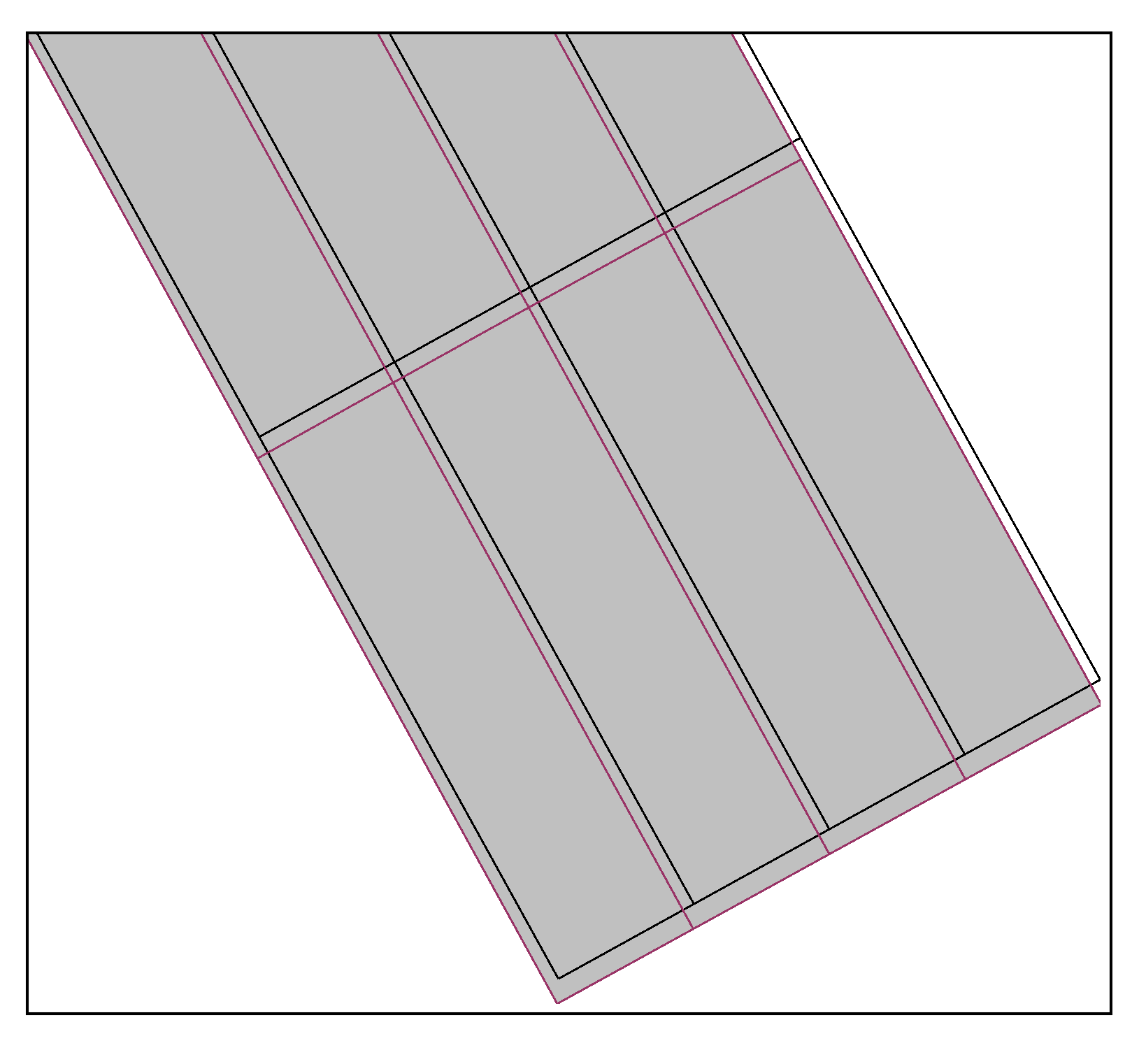} \\
\end{array}$
\end{center}
% \vspace{-7mm}
\caption{Displacement comparison of the initial and optimized quadratic meshes with the reference solution for the 2D cantilever beam. The adapted mesh (shown in purple) captures the physical response (reference solution in solid gray) with higher fidelity than the initial mesh (shown in black).}
\label{fig_beam_morphed_p2}
\end{figure}

The domain $\Omega$ for the elastic beam example has a length $L=1.0$ and height $H=0.1$, and is discretized with $10 \times 4$ uniform-sized Cartesian-aligned quadrilateral elements. The linear elastic isotropic material is characterized by Young's modulus of $E=1.0$ and Poisson ratio of $\nu = 0.3$. A zero Dirichlet boundary condition constrains the left face of $\Omega$ while the top face of $\Omega$ is subject to a pressure load of ${t}=-1.0 \boldsymbol{e}_2$. To minimize the energy norm of the discrete error, we minimize the negative compliance as defined in Section~\ref{sec::compliance}. We choose weight $\alpha = 10^{6}$, and we compare the maximal tip displacement computed with the initial and optimized mesh to the reference displacement in Table~\ref{tab:mesh_comparison_beam}. The reference solution is obtained in an enriched space using fourth-order finite element bases on a mesh that is refined twice compared to the input mesh. The warped displacement field is visualized for a linear $p=1$ and quadratic $p=2$ mesh in Figures \ref{fig_beam_disp_mag_p1} and \ref{fig_beam_morphed_p2}, respectively. For the linear case, we observe an improvement of more than $ 30 \% $.

The optimized mesh results in an improved load functional, $l(\vec{u}_h)$, for both quadratic and linear shape functions, closer to that of the reference solution in comparison to the corresponding initial meshes. The improvement is less significant for quadratic shape functions as they are more capable of representing the displacement of the beam. Observations of improvements are similarly made for both the tip $\vec{u}$ and average displacements $\mid \vec{u} \mid$. These results convey the capability of the mesh $r$-adaptivity framework to produce meshes with a smaller finite element error.

\begin{table}[htbp]
    \centering
    \begin{tabular}{|l|c|c|c|}
        \hline
        \textbf{Mesh} & $\vec{u}$ & $\mid \vec{u} \mid$ & $l\left( u_h \right)$ \\
        \hline \hline
        Reference & (-0.0907,-1.3766) &  1.3796 & $5.55\times10^{-4}$\\
        \hline \hline
        Initial (p=1) & (-0.0662,-0.9999) & 1.0021 & $4.02\times10^{-4}$ \\
        \hline
        Optimized (p=1) & (-0.0736,-1.1776) & 1.1799 & $4.89\times10^{-4}$\\
        \hline
        \hline
        Initial (p=2) & (-0.0903,-1.3679) & 1.3708 & $5.50\times10^{-4}$ \\
        \hline
        Optimized (p=2) & (-0.0907,-1.3759) & 1.3789 &  $5.54\times10^{-4}$ \\
        \hline
    \end{tabular}
    \caption{Comparison of maximum displacement ($\vec{u}$ and $\mid \vec{u} \mid$) at the bottom right corner and load functional ($l(\vec{u}_h)$) for the initial mesh and optimized mesh with the reference solution.}
    \label{tab:mesh_comparison_beam}
\end{table}

\begin{figure}[htbp]
\begin{center}
$\begin{array}{ccc}
\includegraphics[width=0.1\textwidth]
{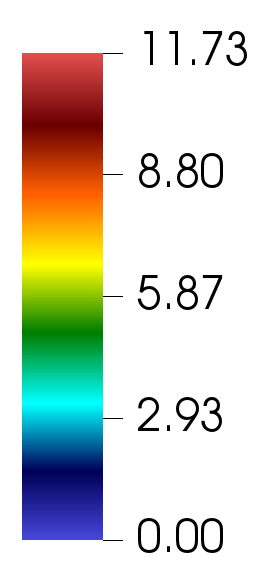} &
\includegraphics[width=0.4\textwidth]
{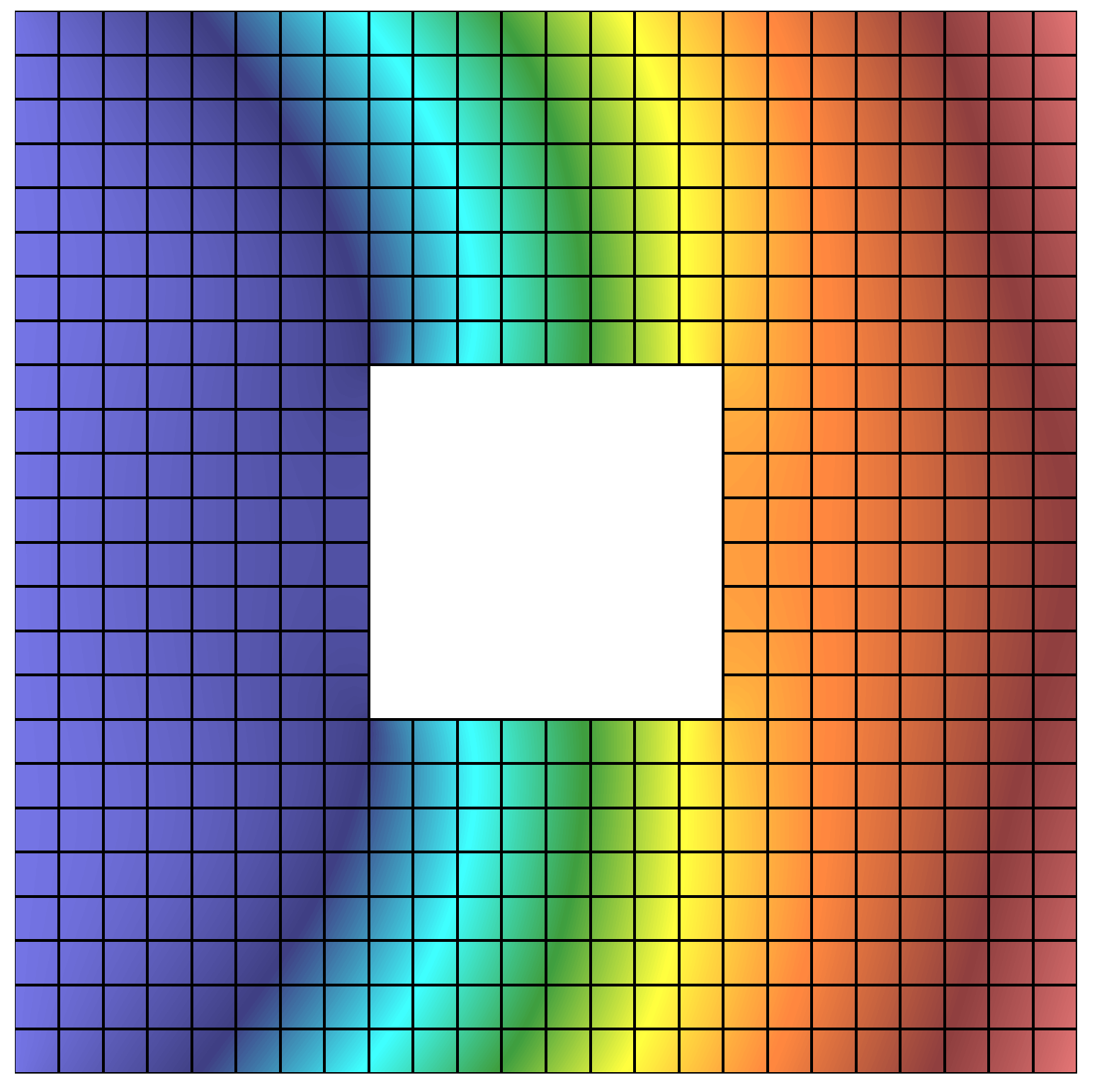} &
\includegraphics[width=0.4\textwidth]
{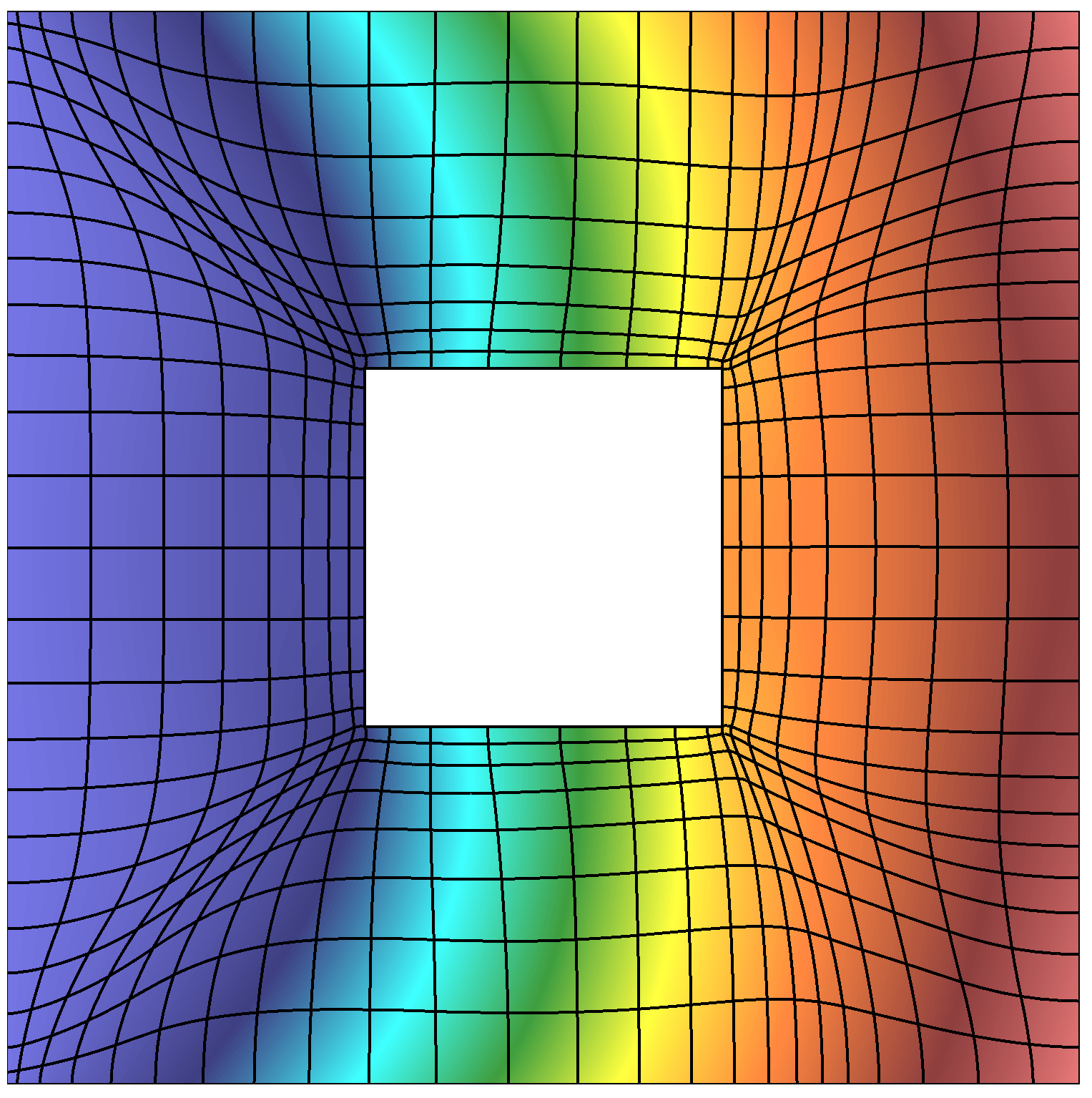} \\
\end{array}$
\end{center}
% \vspace{-7mm}
\caption{Initial and optimized mesh for the shear wall example.}
\label{fig_plate_hole_morphed_p2}
\end{figure}

Figure \ref{fig_plate_hole_morphed_p2} shows the initial mesh for the shear wall example, which is discretized with 512 elements in the domain defined as $\Omega=[0,1]^2 \setminus [1/3,2/3]^2$. The optimization process deforms the elements in the shear band regions formed around the re-entrant corners of the square hole. The displacement field is relatively uniform and aligned with the elements in the rest of the computational domain so that the elements preserve their rectangular shape. The initial and optimized mesh solutions are obtained with second-order elements, and the reference solution is obtained on a quartic mesh with 131,072 elements. A comparison between the maximum displacement is shown in Table \ref{tab:mesh_comparison}. As expected, The optimized mesh results in an improved load functional closer to that of the reference solution. Similar observations of improvement are made for the tip displacement $\vec{u}$ and the average displacement $\mid \vec{u} \mid$.

\begin{table}[htbp]
    \centering
    \begin{tabular}{|l|c|c|c|}
        \hline
        \textbf{Mesh} & $\vec{u}$ & $\mid \vec{u} \mid$ & $l\left(u_h \right)$ \\
        \hline \hline
        Reference & (-3.60617,-11.1919) & 11.7585 & 10.6256 \\
        \hline \hline
        Initial & (-3.58455,-11.0965) & 11.6611 & 10.5464\\
        \hline
        Optimized & (-3.59557,-11.1684) & 11.7329 & 10.6124\\
        \hline
    \end{tabular}
    \caption{Comparison of maximum displacement ($\vec{u}$ and $\mid \vec{u} \mid$) at the bottom right corner and load functional ($l(\vec{u}_h)$) for the initial mesh and optimized mesh with the reference solution.}
    \label{tab:mesh_comparison}
\end{table}

\section{Conclusions}
\label{sec::Conclusion}

In this work, we introduced a general and flexible framework for PDE-constrained $r$-adaptivity of high-order meshes in the context of finite element discretizations. By formulating mesh optimization as a variational problem that couples mesh quality and PDE solution error measure, our method enables control over both discretization error and mesh quality. The approach leverages adjoint sensitivity analysis for efficient gradient computation and employs the MMA solver to robustly solve the resulting large-scale optimization problem. Numerical examples involving elliptic PDEs, including the Poisson and linear elasto-static equations, demonstrate the method’s applicability to different element types, spatial dimensions, and objective functions.

The results show that our method effectively concentrates mesh resolution in regions most critical to solution accuracy as defined by the prescribed error measures.
The use of TMOP ensures high mesh quality throughout the optimization process, and convolution-based gradient regularization improves robustness.

Future work will incorporate this method into differentiable physical simulations and automated design optimization workflows. Additional research directions include extending the methodology to time-dependent PDEs, integrating advanced error estimators, and evaluating its performance in large-scale, GPU-accelerated environments.

\section*{Acknowledgments} The authors would like to thank Veselin Dobrev for his input on the load functional-based error measure in Section \ref{sec::compliance}. This work was performed under the auspices of the U.S. Department of Energy by Lawrence Livermore National Laboratory under Contract DE-AC52-07NA27344. LLNL-JRNL-2006068.

\bibliography{lit}

\end{document}